\documentclass[a4paper,12pt,reqno,twosided]{amsart} 
\usepackage{amsmath,amssymb,amsthm}
\usepackage{graphpap,latexsym,epsf}
\usepackage{epsf}
\usepackage{color,graphics}
\usepackage[numbers,sort&compress,longnamesfirst]{natbib}

\theoremstyle{plain}
\newtheorem{theorem}{Theorem}[section]  
\newtheorem{corollary}[theorem]{Corollary}  
\newtheorem{lemma}[theorem]{Lemma}  
\newtheorem{proposition}[theorem]{Proposition}
\newtheorem*{theorem*}{Theorem}

\theoremstyle{plain}  
\newtheorem{remark}[theorem]{Remark}
\newtheorem{definition}[theorem]{Definition}

\newenvironment{pfn}[1]{\noindent{\bf Proof of {#1}\enspace}}{
\hfill$\Box$\medskip}

\newtheoremstyle{citing}
  {3pt}
  {3pt}
  {\itshape}
  {}
  {\bfseries}
  {.}
  {.5em}
  {\thmnote{#3}}

\theoremstyle{citing}

\numberwithin{equation}{section}

\newcommand{\dist}{\text{dist}}

\newcommand{\rand}{\partial} 
\newcommand{\where}{\,|\,}

\newcommand{\sgn}{\text{\rm sgn}}

\newcommand{\intg}{\mathop{\int}\limits}

\newcommand{\laplace}{\Delta}

\newcommand{\nz}{{\mathbb N}}

\newcommand{\rz}{{\mathbb R}}

\newcommand{\eps}{\varepsilon}  
\renewcommand{\phi}{\varphi}

\renewcommand{\div }{{\rm div}\,}
\renewcommand{\a }{\alpha } 
\renewcommand{\b }{\beta } 
\renewcommand{\d}{\delta }
\newcommand{\smallo}{{o(1)}}
\newcommand{\ui }{u_i} 
\newcommand{\media}{\mkern12mu\hbox{\vrule height4pt depth-3.2pt
    width5pt} \mkern-16mu\int}
\newcommand{\mediapiccola}{\mkern5mu\hbox{\vrule height4pt depth-3.5pt width4pt} \mkern-12.5mu\int}

\newcommand{\Di }{\mathcal{D}^{1,2}_{a}(\R^N) } 
\newcommand{\e }{\varepsilon }

\renewcommand{\l }{\lambda }

\newcommand{\n }{\nabla }

\newcommand{\R}{\mathbb{R}}

\begin{document}
 
\title[Degenerate critical elliptic equations]
{Compactness and existence results for degenerate critical elliptic equations}

\author{Veronica Felli}
\author{Matthias Schneider}
\thanks{ V. F. is supported by M.U.R.S.T. under the national project
  ``Variational Methods and Nonlinear Differential Equations''. During
  the preparation of the paper M. S. research was supported by a
  S.I.S.S.A. postdoctoral fellowship.} 
\address{Scuola Internazionale Superiore di Studi Avanzati\\
S.I.S.S.A.\\
Via Beirut 2-4\\
34014 Trieste, Italy}
\email{felli@sissa.it} 
\address{Mathematisches Institut\\
    Im Neuenheimer Feld 288\\
   D-69120 Heidelberg, Germany}
\email{mschneid@mathi.uni-heidelberg.de} 

\date{June 25, 2003}  
\keywords{critical exponents, blow-up analysis, Caffarelli-Kohn-Nirenberg inequality, Leray-Schauder
  degree}
\subjclass{35J70, 35B33, 35B45}

\begin{abstract}
This paper is devoted to the study of degenerate critical
elliptic equations of Caffarelli-Kohn-Nirenberg type. By means of
blow-up analysis techniques, we prove an a-priori estimate in a weighted
space of continuous functions. From this compactness result, the
existence of a solution to our problem is proved by
exploiting the homotopy invariance of the Leray-Schauder degree. 
\end{abstract}

\maketitle

\section{Introduction}
\label{s:intro}
\noindent
We will consider the following equation in $\R^N$ in dimension
$N\geq3$, which is a prototype of more general nonlinear degenerate elliptic equations
describing anisotropic physical phenomena,
\begin{align}
\label{eq:45}
-\div (|x|^{-2\a}\n v)-\displaystyle{\frac{\l}{|x|^{2(1+\a)}}}\,v=
 K(x)\frac{v^{p-1}}{|x|^{\b p}},\quad v\ge 0,\quad v \in
 {\mathcal{D}^{1,2}_{\a}(\R^N)}\setminus \{0\}, 
\end{align}
where $K \in C^2(\rz^N) \cap L^\infty(\rz^N)$ is positive and
\begin{align}
\label{eq:46}
&\a<\frac{N-2}{2},\quad \a\le \b< \a+1,\\
\label{eq:53}
&\l<\bigg(\frac{N-2-2\a}2\bigg)^2,\quad
p=p(\a,\b)=\displaystyle{\frac{2N}{N-2(1+\a-\b)}}.
\end{align}
We look for weak solutions in ${\mathcal{D}^{1,2}_{\a}(\R^N)}$ defined as the completion
of $C^{\infty}_0(\R^N)$ with respect to the norm
\[
\|u\|_{{\mathcal{D}^{1,2}_{\a}(\R^N)}}:=\left[\int_{\R^N}|x|^{-2\a}|\nabla
  u|^2\,dx\right]^{1/2}.
\] 
The range of $\a$, $\b$ and the definition of $p$ are
related to Caffarelli-Kohn-Nirenberg 
inequalities, denoted by CKN-inequalities in the sequel, (see
\cite{CKN,CatrinaWang} and the references therein), as 
for any $\a$, $\b$ satisfying (\ref{eq:46}) there exists exactly one exponent $p=p(\a,\b)$
such that  
\begin{align}
\label{eq:CKN}
\left(\int_{\R^N}|x|^{-\b p}|u|^p\,dx\right)^{2/p}\leq {\mathcal
  C}_{\a,\b}\int_{\R^N}|x|^{-2\a}|\n u|^2\,dx \qquad \forall u \in C_0^\infty(\rz^N).
\end{align}
Since we are looking for nontrivial nonnegative solutions we must necessarily have that
the quadratic form 
$$Q(\varphi,\varphi):= \int_{\rz^N} |x|^{-2\a}|\nabla \varphi|^2-\lambda
|x|^{-2(1+\a)}|\varphi|^2$$
is positive, that is $\lambda$ has to be smaller than $(N-2-2\a)^2/4$ the best constant in the
related Hardy-type CKN-inequality for $\b=\a+1$ and $p=2$. 
Let us define
\begin{align}\label{eq:ab}
a(\a,\l):=\frac{N-2}2-{\sqrt{\left(\frac{N-2-2\a}2\right)^2-\l}}\text{
  and }
b(\a,\b,\l):=\b+a(\a,\l)-\a. 
\end{align}
The change of variable $u(x)=|x|^{a-\a}v(x)$ shows that equation
(\ref{eq:45}) is equivalent 
to 
\begin{equation}\label{eq:49}
-\div(|x|^{-2a}\n u)=K(x)\frac{u^{p-1}}{|x|^{bp}}, \quad u\ge 0,\quad u \in \Di\setminus\{0\},
\end{equation}
where $a=a(\a,\l)$ and $b=b(\a,\b,\l)$, see Lemma \ref{l:change} of
the Appendix. 
Clearly, if we replace $\a$ by $a$ and $\b$ by $b$ then (\ref{eq:46})-(\ref{eq:53}) still
hold and $p(\a,\b)=p(a,b)$. We will write in the sequel for short that $a$, $b$ and $p$
satisfy (\ref{eq:46})-(\ref{eq:53}).
We will mainly deal with equation (\ref{eq:49}) and look for weak solutions in 
$\Di$. The advantage of working with (\ref{eq:49}) instead of (\ref{eq:45}) is that we know from \cite{SF_reg} 
that weak solutions of (\ref{eq:49}) are H\"older-continuous in $\rz^N$ whereas
solutions to (\ref{eq:45}), as our analysis shows, behave (possibly singular) like $|x|^{\a-a}$ at the origin.  
The main difficulty in facing problem (\ref{eq:49}) is the lack of compactness as $p$ is the
critical exponent in the related CKN-inequality.
More precisely, if $K$ is a positive constant equation (\ref{eq:49}) is invariant under the action
of the non-compact group of dilations, in the sense that 
if $u$ is a solution of (\ref{eq:49}) then for any positive $\mu$ the dilated 
function
\[
\mu^{-\frac{N-2-2a}{2}}u(x/\mu)
\] 
is also a solution with the same norm in $\Di$. The dilation invariance, as we see in
(\ref{eq:59}) below, gives rise to a
non-compact, one dimensional manifold of solutions for $K \equiv K(0)$.  
\par
Our first theorem provides sufficient conditions on $K$ ensuring
compactness of the set of solutions by means of an a-priori bound in a weighted space $E$
defined by 
$$
E:=\Di\cap C^0(\R^N,(1+|x|^{N-2-2a})),
$$
where 
$$
C^0(\R^N,(1+|x|^{N-2-2a})):=\{u\in C^0(\R^N):\, u(x)(1+|x|^{N-2-2a})\in
L^{\infty}(\R^N)\}
$$
is equipped with the norm
\[
\|u\|_{C^0(\R^N,(1+|x|^{N-2-2a}))}:=\sup_{x\in \R^N} |u(x)|(1+|x|^{N-2-2a}).
\]
We endow $E$ with the norm
\[
\|u\|_E=\|u\|_{\Di}+\|u\|_{C^0(\R^N,(1+|x|^{N-2-2a}))}.
\]
The uniform bound in $E$ of the set of solutions to (\ref{eq:49}) will provide the necessary
compactness needed in the sequel. We formulate the compactness result in terms of $\a$, $\beta$ and
$v$ the parameters of equation (\ref{eq:45}), where we started from.
Let us set 
\begin{equation}\label{eq:68}
   \tilde K(x):=K(x/|x|^2).
\end{equation}
\begin{theorem}{\bf(Compactness)} \label{t:compactness}
Let $\a,\b,\l$ satisfy (\ref{eq:46})-(\ref{eq:53}) and
\begin{align}
\label{eq:47}
&\l\ge -\a(N-2-\a),\\
\label{eq:58}
&\bigg(\frac{N-2-2\a}{2}\bigg)^2-1<\l,\\
\label{eq:48}
&\b>\a,\quad p>\frac{2}{\sqrt{\left(\frac{N-2-2\a}2\right)^2-\l}}.
\end{align}
Suppose $K\in C^2(\R^N)$ satisfies 
\begin{align}
\label{eq:60}
&\tilde K \in C^2(\rz^N), \text{ where } \tilde K(x)\text{ is defined
  in (\ref{eq:68})},\\
\label{eq:7}
&\n K(0)=0,\quad \Delta K(0)\not=0,\quad\text{and}\quad 
\n \tilde K(0)=0,\quad \Delta \tilde K(0)\not=0,
\end{align}
and for some positive constant $A_1$
\begin{align}
\label{eq:51}
&1/A_1\leq K(x),\quad\forall\, x\in\R^N.
\end{align}
Then there is $C_K>0$ such that for any $t\in (0,1]$ and any solution $v_t$ of 
\addtocounter{equation}{1}
\begin{align}
\label{eq:52}
\tag*{$(\theequation)_{t}$}
\begin{split}
-\div (|x|^{-2\a}\n v)-\displaystyle{\frac{\l}{|x|^{2(1+\a)}}}\,v
=(1+t\,(K(x)-1))\frac{v^{p-1}}{|x|^{\b p}},\\ 
\ v\ge0,\ v \in {\mathcal{D}^{1,2}_{\a}(\R^N)}\setminus \{0\},
\end{split}  
\end{align}
we have $\||x|^{a-\a}v_t\|_E< C_K$ and 
\begin{align}
\label{eq:50}
C_K^{-1}<|x|^{a-\a}(1+|x|^{N-2-2a})v_t(x)<C_K \text{ in } \R^N\setminus\{0\}.
\end{align}
\end{theorem}
To prove the above compactness result we adapt the arguments of \cite{li95} to carry out a fine
blow-up analysis for (\ref{eq:49}). Assumptions (\ref{eq:47})-(\ref{eq:48}) imply
\[
\begin{array}{lll}
&(\ref{eq:47}) \ \Longrightarrow\   a\ge 0,  &(\ref{eq:58})\  \Longrightarrow\ 
\frac{N-4}{2}<a<\frac{N-2}{2} \\[5pt]
&(\ref{eq:48})\  \Longrightarrow\ \frac{4}{N-2-2a}< p<2^*=\frac{2N}{N-2}.
\end{array}
\]
A key ingredient is the exact knowledge of the
solutions to the limit problem with $K \equiv {\rm const}$, which is only available for
$a\ge 0$. In \cite{ChouChu} (see also
\cite{terracini}) it is shown through 
the method of moving planes that if $a\ge 0$ then any locally bounded positive solution in
$C^2(\R^N\setminus\{0\})$ of~(\ref{eq:49}) with $K \equiv K(0)$ is of the form
\begin{align}
\label{eq:59}
z_{K(0),\mu}^{a,b}:=\mu^{-\frac{N-2-2a}2}z_{K(0)}^{a,b}\Big(\frac
x{\mu}\Big),\quad\mu>0,   
\end{align}
where $z_{K(0)}^{a,b}=z_1^{a,b}\big(x\, K(0)^{\frac 2{(p-2)(N-2-2a)}}\big)$
and $z_1^{a,b}$ is explicitly given by  
\begin{align*}
&z_1^{a,b}(x)=
\left[1+\frac{N-2(1+a-b)}{N(N-2-2a)^2} |x|^{\frac{2(1+a-b)(N-2-2a)}{N-2(1+a-b)}}\right]^{-\frac{N-2(1+a-b)}{2(1+a-b)}}.
\end{align*}
For $a<0$ the set of positive solutions becomes more and more complicated as $a\! \to\!
-\infty$ due to the existence of non-radially symmetric solutions (see
\cite{FS1,CatrinaWang,CatrWang01}).  
Up to now, our blow-up analysis is
only available for~$p<2^*$; the case $p=2^*$ presents additional
difficulties because besides the blow-up profile~$z_1^{a,b}$ a second blow-up profile described 
by the usual Aubin-Talenti instanton of Yamabe-type equations may occur. 
The further restrictions on $a$, $p$ and $K$ should be compared to
the so-called flatness-assumptions in problems of prescribing scalar curvature.\par
 Non-existence
results for equation (\ref{eq:49}) can be obtained using a
Pohozaev-type identity, i.e. any solution $u$ to \eqref{eq:49} satisfies the following identity
\[
\int_{\R^N}(\n K(x)\cdot x)\frac{u^p}{|x|^{bp}}\,dx=0,
\]
provided the integral
is convergent and $K$ is bounded and smooth enough (see Corollary
\ref{c:pohoRn}). This implies that there are no such solutions if
$\n K(x)\cdot x$ does not change
sign in $\R^N$ and $K$ is not constant.\par 
The above compactness result allows us to exploit the homotopy
invariance of the Leray-Schauder degree to pass from $t$ small 
to $t=1$ in~\ref{eq:52}. We compute the degree of positive solutions
to \ref{eq:52} for small $t$ using a Melnikov-type function introduced
in \cite{AmBa1,AmBa2} and show that it equals (see Theorem \ref{t:grado})
\[
-\frac{\sgn\Delta K(0)+\sgn\Delta \tilde K(0)}2.
\]  
In particular, we prove the following existence result.    
\begin{theorem}{\bf(Existence)} \label{t:existence}
Under the assumptions of Theorem \ref{t:compactness}, if, moreover,  $p>3$ and
\[
\sgn\Delta K(0)+\sgn  \Delta \tilde K(0) \neq 0
\]
then equation (\ref{eq:45}) has a positive solution $v$ such that
$|x|^{a-\a}v\in B_{C_K}(0)\subset E$ and $v$ satisfies (\ref{eq:50}). 
\end{theorem}
\noindent 
The assumption $p>3$ is essentially technical and yields $C^3$
regularity of the functional associated to the problem which is needed
in the computation of the degree.\par
In \cite{FS1} problem (\ref{eq:45}) is studied in the case in which
$K$ is a small perturbation of a constant, i.e. in the case $K=1+\e
k$, using a perturbative method introduced in~\cite{AmBa1,AmBa2}. We extend some of the results
in \cite{FS1} to the nonperturbative case.
Problem~(\ref{eq:45}) for~$\a=\b=0$ (hence $p=2^*$) and
$0<\l<(N-2)^2/4$ is treated by Smets 
\cite{smets} who proves that in dimension $N=4$ there exists a
positive solution provided $K\in C^2$ is positive and
$K(0)=\lim_{|x|\to\infty}K(x)$. Among other existence and multiplicity
results, in \cite{AbdelFelliPeral} positive solutions to (\ref{eq:45})
for $\a=\b=0$, $p=2^*$, and $0<\l<(N-2)^2/4$ 
are found via the concentration
compactness argument, under assumptions ensuring that the
mountain-pass level stays below the 
compactness threshold at which Palais-Smale condition fails. 
We emphasize that the solution we find in
Theorem \ref{t:existence} can stay above such a threshold.\par
\begin{remark}\label{rem:intro2} 
If we drop the assumption $\alpha<\frac{N-2}2$ we may still change the variables 
$u(x)=|x|^{a-\a}v(x)$, where $a$ is given in (\ref{eq:ab}), 
and we still obtain weak solutions $u$ of (\ref{eq:49}) in ${\mathcal{D}^{1,2}_{a}(\R^N)}$.
But in this case the transformation $v(x)=|x|^{\a-a}u(x)$ gives rise only to
classical solutions of~(\ref{eq:45}) in $\R^N\setminus\{0\}$ but not to 
distributional solutions in the whole $\R^N$.
\end{remark}
The paper is organized as follows. In Section 2 we prove a Pohozaev
type identity for equation (\ref{eq:49}). In Section 3 we introduce
the notion of isolated and isolated simple blow-up point which was
first introduced by Schoen \cite{schoen91} and provide the main local
blow-up analysis. In Section 4 we prove Theorem \ref{t:compactness} by
combining the Pohozaev type identity with the results of our local
blow-up analysis. Section 5 is devoted to the computation of the
Leray-Schauder degree and to the proof of the existence
theorem. Finally in the Appendix we collect some technical lemmas.

\begin{center}

{\bf Acknowledgements}

\end{center}

The authors would like to thank Prof. A. Ambrosetti for his interest
in their work and for helpful suggestions.


\section{A Pohozaev-type identity}
\begin{theorem}\label{t:pohozaev}
Let $\Omega$ be a bounded domain in $\R^N$ ($N\geq 3$) with smooth
boundary, $a$, $b$, and $p$ satisfy (\ref{eq:46})-(\ref{eq:53}), $K\in
C^1(\overline\Omega)$ and $u \in \Di$ be a weak
 positive solution of  
\begin{equation}\label{eq:poho1}
-\div(|x|^{-2a}\n u)= K(x)\frac{u^{p-1}}{|x|^{bp}}, \quad x\in\Omega.
\end{equation}
There holds
\begin{align*}
&\frac {1}p \int_{\Omega}(\n K\cdot x)\,\frac{u^p}{|x|^{bp}}\,
dx-\frac {1}p\int_{\partial \Omega}K(x)\frac{u^p}{|x|^{bp}}  x\cdot
  \nu=\frac{N-2-2a}2\int_{\partial \Omega}|x|^{-2a}u\n u\cdot
  \nu\\
&\quad-\frac 12 \int_{\partial \Omega}|x|^{-2a}|\n
  u|^2x\cdot\nu +\int_{\partial \Omega}|x|^{-2a}(x\cdot\n u)(\n u\cdot\nu)
\end{align*}
where $\nu$ denotes the unit normal of the boundary.
\end{theorem}
\begin{proof}
Note that 
\[
\int_{0}^1ds\int_{\partial
  B_s(0)}\bigg[\frac{|K(x)|u^p}{|x|^{bp}}+\frac{|\nabla
  u|^2}{|x|^{2a}}\bigg]=\int_{B_1(0)}\bigg[\frac{|K(x)|u^p}{|x|^{bp}}
+\frac{|\nabla u|^2}{|x|^{2a}}\bigg]<\infty  
\]
which implies that there exists a sequence $\e_n\to 0^+$ such that 
\begin{equation}\label{eq:poho7}
\e_n\int_{\partial
  B_{\e_n}(0)}\bigg[\frac{|K(x)|u^p}{|x|^{bp}}+\frac{|\nabla
  u|^2}{|x|^{2a}}\bigg]\to 0 
\end{equation}
as $n\to \infty$. Let $\Omega_{\e_n}:=\Omega\setminus B_{\e_n}(0)$. 
Multiplying equation (\ref{eq:poho1}) by $x\cdot\n u$ and integrating
over $\Omega_{\e_n}$ we obtain
\begin{equation}\label{eq:poho2}
-\sum_{j,k=1}^N\int_{\Omega_{\e_n}}\frac{\partial }{\partial
 x_j}\bigg(|x|^{-2a}\frac{\partial u}{\partial
 x_j}\bigg)x_k\,\frac{\partial u}{\partial x_k}\,
 dx= \sum_{k=1}^N\int_{\Omega_{\e_n}}x_k\,\frac{\partial u}{\partial x_k}\,
 K(x)\frac{u^{p-1}}{|x|^{bp}}\,dx.
\end{equation}
Let us first consider the right-hand side of (\ref{eq:poho2}). Integrating
by parts we have
\begin{align}\label{eq:poho3}
&\sum_{k=1}^N\int_{\Omega_{\e_n}}x_k\,\frac{\partial u}{\partial x_k}\,
 K(x)\frac{u^{p-1}}{|x|^{bp}}\,dx=\bigg(b-\frac
 Np\bigg)\int_{\Omega_{\e_n}}K(x)\frac{u^p}{|x|^{bp}}\,dx \notag\\
&\qquad-\frac
 1p\sum_{k=1}^N\int_{\Omega_{\e_n}}u^px_k\frac{\partial K}{\partial
 x_k}|x|^{-bp}+\frac
 1p\sum_{k=1}^N\int_{\partial
 \Omega_{\e_n}}x_k\nu_k\frac{K(x)u^p}{|x|^{bp}}.
\end{align}
Integrating by parts in the left-hand side of (\ref{eq:poho2}), we obtain
\begin{align}\label{eq:poho4}
&-\sum_{j,k=1}^N\int_{\Omega_{\e_n}}\frac{\partial }{\partial
 x_j}\bigg(|x|^{-2a}\frac{\partial u}{\partial
 x_j}\bigg)x_k\,\frac{\partial u}{\partial x_k}\,
 dx =-\frac{N-2-2a}2\int_{\Omega_{\e_n}}|x|^{-2a}|\n u|^2\,dx\notag\\
&\qquad +\frac 12
 \int_{\partial \Omega_{\e_n}}|x|^{-2a}|\n u|^2 x\cdot\nu
 -\int_{\partial\Omega_{\e_n}}|x|^{-2a}(x\cdot\n u)(\n u\cdot 
 \nu).
\end{align}
From (\ref{eq:poho2}), (\ref{eq:poho3}), and (\ref{eq:poho4}), we have
\begin{align*}
& \bigg(b-\frac Np\bigg)\int_{\Omega_{\e_n}}K(x)\frac{u^p}{|x|^{bp}}\,dx- \frac
 1p\sum_{k=1}^N\int_{\Omega_{\e_n}}u^px_k\frac{\partial K}{\partial
 x_k}|x|^{-bp}\notag\\
&\qquad\quad + \frac
 1p\sum_{k=1}^N\int_{\partial
 \Omega_{\e_n}}x_k\nu_k\frac{K(x)u^p}{|x|^{bp}}\\
&\qquad =-\frac{N-2-2a}2\int_{\Omega_{\e_n}}|x|^{-2a}|\n u|^2\,dx+\frac 12
 \int_{\partial \Omega_{\e_n}}|x|^{-2a}|\n u|^2 x\cdot\nu\notag\\
&\qquad\quad -\int_{\partial\Omega_{\e_n}}|x|^{-2a}(x\cdot\n u)(\n u\cdot
 \nu).
\end{align*}
Because of the integrability of $|x|^{-bp}u^p$ and of $|x|^{-2a}|\n u|^2$, it is clear that 
\begin{align*}
&\bigg(b-\frac Np\bigg)\int_{\Omega_{\e_n}}K(x)\frac{u^p}{|x|^{bp}}\,dx-\frac
 1p\sum_{k=1}^N\int_{\Omega_{\e_n}}u^px_k\frac{\partial K}{\partial
 x_k}|x|^{-bp}\\
&\qquad\mathop{\longrightarrow}\limits_{\e\to 0^+} \bigg(b-\frac Np\bigg)\int_{\Omega}K(x)\frac{u^p}{|x|^{bp}}\,dx-\frac
 1p\sum_{k=1}^N\int_{\Omega}u^px_k\frac{\partial K}{\partial
 x_k}|x|^{-bp}
\end{align*}
and 
\[
\int_{\Omega_{\e_n}}|x|^{-2a}|\n u|^2\,dx\mathop{\longrightarrow}\limits_{\e\to 0^+}\int_{\Omega}|x|^{-2a}|\n u|^2\,dx.
\]
Hence, in view of (\ref{eq:poho7}), we have
\begin{align}\label{eq:poho8}
& \bigg(b-\frac Np\bigg)\int_{\Omega}K(x)\frac{u^p}{|x|^{bp}}\,dx- \frac
 1p\sum_{k=1}^N\int_{\Omega}u^px_k\frac{\partial K}{\partial
 x_k}|x|^{-bp}\notag\\
&\qquad\quad + \frac
 1p\sum_{k=1}^N\int_{\partial
 \Omega}x_k\nu_k\frac{K(x)u^p}{|x|^{bp}}\notag\\
&\qquad =-\frac{N-2-2a}2\int_{\Omega}|x|^{-2a}|\n u|^2\,dx+\frac 12
 \int_{\partial \Omega}|x|^{-2a}|\n u|^2 x\cdot\nu\notag\\
&\qquad\quad -\int_{\partial\Omega}|x|^{-2a}(x\cdot\n u)(\n u\cdot
 \nu).
\end{align}
Multiplying equation (\ref{eq:poho1}) by $u$ and integrating by parts, we
have
\begin{equation}\label{eq:poho9}
\int_{\Omega}|x|^{-2a}|\n u|^2\,dx=\int_{\partial
  \Omega}|x|^{-2a}u\,\frac{\partial u}{\partial
  \nu}+ \int_{\Omega}K(x)\frac{u^p}{|x|^{bp}}\,dx.
\end{equation}
The conclusion follows from (\ref{eq:poho8}), (\ref{eq:poho9}), and
from the identity $\frac{N-bp}{p}-\frac{N-2-2a}2=0$.
\end{proof}

\begin{corollary}\label{c:poho1}
If $a$, $b$, and $p$ satisfy (\ref{eq:46})-(\ref{eq:53}),  $K\in
C^1(\overline B_{\sigma})$ and $u$ be a weak
 positive solution in $\Di$ of  
\begin{equation}\label{eq:poho6}
-\div(|x|^{-2a}\n u)= \,K(x)\frac{u^{p-1}}{|x|^{bp}}, \quad x\in
 B_{\sigma}:=\{x\in\R^N:\ |x|<\sigma\}
\end{equation}
then
\begin{align}\label{eq:29}
\frac {1}p \int_{B_{\sigma}}(\n K\cdot x)\,\frac{u^p}{|x|^{bp}}\,
dx- \,\frac {\sigma}p\int_{\partial
  B_{\sigma}}K(x)\frac{u^p}{|x|^{bp}} =\int_{\partial
  B_{\sigma}}B(\sigma, x, u, \n u) 
\end{align}
where
\[
B(\sigma, x, u, \n u)=\frac{N-2-2a}2\,|x|^{-2a}u\,\frac{\partial
  u}{\partial \nu}-\frac {\sigma}2\, |x|^{-2a}|\n
  u|^2+\sigma \,|x|^{-2a}\bigg(\frac{\partial u}{\partial \nu}\bigg)^2.
\]
\end{corollary}
\begin{corollary}\label{c:pohoRn}
Let $u$ be a weak positive solution in $\Di$ of 
\[
-\div(|x|^{-2a}\n u)= \,K(x)\frac{u^{p-1}}{|x|^{bp}}, \quad x\in\R^N
\]
where $a$, $b$, and $p$ satisfy (\ref{eq:46})-(\ref{eq:53}) and $K\in
L^{\infty}\cap C^1(\R^N)$, $|\n K(x)\cdot x|\leq{\rm 
  const}$. Then 
\begin{equation}
  \label{eq:pohoRn}
\int_{\R^N}(\n K(x)\cdot x)\frac{u^p}{|x|^{bp}}\,dx=0.
\end{equation}
\end{corollary}
\begin{proof}
Since 
\[
\int_{0}^{+\infty}ds\int_{\partial B_s }\bigg[\frac{|K(x)|u^p}{|x|^{bp}}+\frac{|\nabla u|^2}{|x|^{2a}}\bigg]=\int_{\R^N}\bigg[\frac{|K(x)|u^p}{|x|^{bp}}+\frac{|\nabla u|^2}{|x|^{-2a}}\bigg]<\infty
\]
there exists a sequence $R_n\to +\infty$ such that 
\begin{equation}\label{eq:55}
R_n\int_{\partial
  B_{R_n} }\bigg[\frac{|K(x)|u^p}{|x|^{bp}}+\frac{|\nabla
  u|^2}{|x|^{2a}}\bigg]\ \mathop{\longrightarrow}\limits_{n\to\infty}\ 0.
\end{equation}
From Corollary \ref{c:poho1} we have that
\begin{align}\label{eq:56}
\frac {1}p \int_{B_{R_n}}(\n K\cdot x)\,\frac{u^p}{|x|^{bp}}\,
dx=&\frac {R_n}p\int_{\partial B_{R_n}}K(x)\frac{u^p}{|x|^{bp}}
+\frac{N-2-2a}2\int_{\partial B_{R_n}} |x|^{-2a}u\,\frac{\partial
  u}{\partial \nu}\notag\\
&-\frac {R_n}2\int_{\partial B_{R_n}} |x|^{-2a}|\n
  u|^2+R_n \int_{\partial B_{R_n}}|x|^{-2a}\bigg(\frac{\partial
    u}{\partial \nu}\bigg)^2.
\end{align}
In view of (\ref{eq:55}) and noting that from H\"older inequality 
\begin{align*}
\int_{\partial B_{R_n}} |x|^{-2a}u\,\frac{\partial
  u}{\partial \nu}&=R_n^{b-a}\int_{\partial B_{R_n}}
\frac{u}{|x|^b}\cdot\frac{\n u\cdot \nu}{|x|^a}\\
&\leq |{\mathbb
  S}^N|^{\frac{p-2}{2p}}R_n^{b-a+\frac{(N-1)(p-2)}{2p}-\frac 1p-\frac
  12}\left(R_n\int_{\partial
    B_{R_n}}\frac{u^p}{|x|^{bp}}\right)^{\frac 1p}\left(R_n\int_{\partial
    B_{R_n}}\frac{|\n u|^2}{|x|^{2a}}\right)^{\frac 12}\\
&=|{\mathbb
  S}^N|^{\frac{p-2}{2p}}\left(R_n\int_{\partial
    B_{R_n}}\frac{u^p}{|x|^{bp}}\right)^{\frac 1p}\left(R_n\int_{\partial
    B_{R_n}}\frac{|\n u|^2}{|x|^{2a}}\right)^{\frac 12}
\end{align*}
we can pass to the limit in (\ref{eq:56}) thus obtaining the claim.
\end{proof}
\par\noindent
It is easy to check that the boundary term $B(\sigma, x, u, \n u)$ has
the following properties.
\begin{proposition}\label{p:boundterm}
\[
\begin{array}{ll}
(i)\  &\text{For $u(x)=|x|^{2+2a-N}$, $\sigma>0$, $B(\sigma, x, u, \n
  u)=0$ for all $x\in\partial B_{\sigma}$.}\\[5pt]
(ii)\ &\text{For $u(x)=|x|^{2+2a-N}+A+\zeta(x)$, with $A>0$ and $\zeta(x)$
  some function}\\
& \text{differentiable near $0$ satisfying $\zeta(0)=0$, there
  exists $\bar\sigma$ such that}\\[10pt]
&\qquad\qquad B(\sigma, x, u, \n
  u)<0\quad\text{for all}\ x\in\partial B_{\sigma}\text{ and }0<\sigma<\bar\sigma\\[5pt]
&\text{and}\\[5pt]
&\qquad\qquad\displaystyle{\lim_{\sigma\to0}\int_{\partial B_{\sigma}}B(\sigma, x, u, \n
  u)=-\frac{(N-2-2a)^2}2A|{\mathbb S}^{N-1}|}.
\end{array}
\]
\end{proposition}

\section{Local blow-up analysis}


Let $\Omega\subset \R^N$ be a bounded domain, $a$, $b$, and $p$ satisfy (\ref{eq:46})-(\ref{eq:53}), and
$\{K_i\}_i\subset C(\Omega)$ satisfy, for some constant $A_1>0$, 
\begin{equation}\label{eq:blow1}
1/A_1\leq K_i(x)\leq A_1,\quad\forall\, x\in\Omega\quad\text{and}\quad
K_i\to K \text{ uniformly in }\Omega.
\end{equation}
Moreover, we will assume throughout this section that $a \ge 0$.
We are interested in the family of problems 
\begin{align}
\label{eq:pi}
\tag{$P_i$}
\begin{split}
-\div (|x|^{-2a}\n u)= K_i(x)\frac{u^{p-1}}{|x|^{bp}}\quad \text{weakly in}\ \Omega,
\quad u>0 \text{ in}\ \Omega, \ u \in \Di.
\end{split} 
\end{align}
\begin{definition}\label{d:1}
Let $\{\ui\}_i$ be a sequence of solutions of $(P_i)$. We say that $0\in\Omega$
is a blow-up point of $\{\ui\}_i$ if there exists
a sequence $\{x_i\}_i$ converging to $0$ such that
\begin{equation}\label{eq:16}
\ui(x_i)\to+\infty\quad\text{and}\quad
u_i(x_i)^{\frac2{N-2-2a}}|x_i|\to 0\quad\text{as 
}i\to+\infty. 
\end{equation}
\end{definition}
\begin{definition}\label{d:2}
Let $\{\ui\}_i$ be a sequence of solutions of $(P_i)$. The point $0$ is said to
be an isolated blow-up point of $\{\ui\}_i$ if 
there exist $0<\bar r< \dist(0,\partial \Omega)$, $\bar C>0$, and
a sequence $\{x_i\}_i$ converging to $0$ such that
$\ui(x_i)\to+\infty$, $u_i(x_i)^{\frac2{N-2-2a}}|x_i|\to 0$ as $i\to+\infty$,
and for any $x\in B_{\bar
    r}(x_i)$
\[ 
\ui(x)\leq 
\bar C\, |x-x_i|^{-\frac{N-2-2a}2}   
\]
where $B_{\bar r}(x_i):=\{x\in \Omega:\ |x-x_i|<\bar r\}$.
\end{definition}
\noindent If $0$ is an isolated blow-up point of $\{\ui\}_i$ we
define
\[
\bar \ui(r)=\media_{\partial B_r(x_i)}\ui=\frac{1}{|\partial
  B_r(x_i)|}\int_{\partial B_r(x_i)}\ui,\quad r>0
\]
and 
\begin{align}
\label{eq:62}
\bar w_i(r)=
r^{\frac{N-2-2a}2}\bar\ui(r),
\quad r>0.  
\end{align}
\begin{definition}\label{d:4}
The point $0$ is said to
be an isolated simple blow-up point of $\{\ui\}_i$ if it is an
isolated blow-up point and there exist
some positive $\rho\in(0,\bar r)$ independent of~$i$ and $\tilde C>1$ such that
\begin{align}
\label{eq:isolated-simple}
\bar w_i'(r)<0\quad\text{for $r$ satisfying}\ 
\tilde C\,\ui(x_i)^{-\frac
  2{N-2-2a}}\leq r\leq \rho.
\end{align}
\end{definition}
\noindent Let us now introduce the notion of blow-up at infinity. To this aim,
we consider the Kelvin transform,
\begin{equation}\label{eq:kelvin}
\tilde u_i(x)=|x|^{-(N-2-2a)}\ui\bigg(\frac{x}{|x|^2}\bigg),
\end{equation}
which is an isomorphism of $\Di$. 
If $\ui$ solves $(P_i)$ in a neighborhood of $\infty$, i.e. $\Omega=\rz^N\setminus D$ for some compact set $D$,
then $\tilde u_i$ is a solution of $(P_i)$ where $K_i$ is replaced by $\tilde{K_i}(x)=K_i(x/|x|^2)$ and 
$\Omega$ by $\tilde \Omega=\rz^N \setminus \{x/|x|^2 \where x \in D\}$, a neighborhood of $0$.
\begin{definition}\label{d:5}
Let $\{\ui\}_i$ be a sequence of solutions of $(P_i)$ in a
neighborhood of~
$\infty$. We say that $\infty$ is a blow-up point
(respectively an isolated blow-up point, an isolated simple blow-up
point) if $0$ is a blow-up point (respectively an isolated blow-up
point, an isolated simple blow-up point) of the sequence $\{\tilde u_i\}_i$
defined by the Kelvin transform~(\ref{eq:kelvin}).
\end{definition}
\begin{remark}\label{e:infty}
It is easy to see that $\infty$ is a blow-up point of $\{\ui\}_i$ if
and only if there exists a sequence $\{x_i\}_i$ such that
$|x_i|\to\infty$ as $i\to+\infty$ and 
\[
|x_i|^{N-2-2a}\ui(x_i)\ \mathop{\longrightarrow}\limits_{i\to+\infty}\
 \infty\quad\text{and}\quad |x_i| u_i(x_i)^{\frac 2{N-2-2a}}\
 \mathop{\longrightarrow}\limits_{i\to+\infty}\ 0.
\]
\end{remark}
\noindent In the sequel we will use the notation $c$ to denote a positive
constant which may vary from line to line.


\begin{lemma}
\label{l:harnack:zero}
Let $(K_i)_{i \in \nz}$ satisfy (\ref{eq:blow1}), $\{u_i\}_i$ satisfy (\ref{eq:pi}) and
$x_i \to 0$ be an isolated blow up point. Then there is a positive constant $C=C(N,\bar C,A_1)$  
such that for any $0<r<\min(\bar r/3,1)$ there holds
\begin{align}
\label{eq:11}
\max_{x \in B_{2r}(x_i)\setminus B_{r/2}(x_i)}u_i(x) \le C
\min_{x \in B_{2r}(x_i)\setminus B_{r/2}(x_i)}u_i(x).
\end{align}
\end{lemma}
\begin{proof}
We define $v_i(x) := r^{\frac{N-2-2a}{2}} u_i(rx+x_i)$. Then $v_i$ satisfies in $B_3(0)$
\begin{align}\label{eq:57}
0<v_i(x)<\bar C |x|^{-\frac{N-2-2a}{2}},
\end{align}
and
\begin{align*}
-\div\big(|x+r^{-1}x_i|^{-2a} \nabla v_i(x)\big) &= 
-r^{\frac{N-2-2a}{2}+2+2a} \div\big(|\cdot|^{-2a} \nabla u_i(\cdot)\big)(rx+x_i)\\
&= K_i(rx+x_i)~|x+r^{-1}x_i|^{-bp}~ v_i^{p-1}(x),   
\end{align*}
since
\begin{align*}
\frac{N-2-2a}{2}+2+2a-bp-(p-1)\frac{N-2-2a}{2} &= N-p\Big(\frac{N-2(1+a-b)}{2}\Big)= 0. 
\end{align*}
To prove the claim we use a weighted version of Harnack's inequality applied to $v_i$ and
\begin{align*}
-\div\big(|x+r^{-1}x_i|^{-2a} \nabla v_i(x)\big) - W_i(x) v_i(x) =0 \quad\text{ in } B_{9/4}(0)\setminus B_{1/4}(0),
\end{align*}
where $W_i(x):= K_i(rx+x_i)~|x+r^{-1}x_i|^{-bp}~ v_i^{p-2}(x)$. 
From (\ref{eq:57}) the function $v_i$ is uniformly bounded in $B_{9/4}(0)\setminus B_{1/4}(0)$ and 
the claim follows from Harnack's inequality in \cite{Gutierrez89}. We
mention that $|\cdot+r^{-1}x_i|^{-bp}$ belongs to the class of
potentials required in~\cite{Gutierrez89} (see Lemma
\ref{l:harnack_potential} of the Appendix). 
\end{proof}

\begin{proposition}
\label{prop:Bubble:zero}
Let $\{K_i\}_i$ satisfy (\ref{eq:blow1}), $\{u_i\}_i$ satisfy (\ref{eq:pi}) and
$x_i \to 0$ be an isolated blow up point. Then for any $R_i \to \infty, \, \eps_i \to
0^+$, we have, after passing to a subsequence 
that:\\
\begin{align}
\label{eq:6}
R_i u_i(x_i)^{-\frac{2}{N-2-2a}} \to 0 \text{ as } i \to \infty,  
\end{align}
\begin{align}
\label{eq:4}
\|u_i(x_i)^{-1} u_i\big(u_i(x_i)^{-\frac{2}{N-2-2a}}\cdot + x_i\big) 
- z_{K(0)}^{a,b}(\cdot)\|_{C^{0,\gamma}(B_{2R_i}(0))} \le \eps_i,  
\end{align}
\begin{align}
\label{eq:15}
\|u_i(x_i)^{-1} u_i\big(u_i(x_i)^{-\frac{2}{N-2-2a}}\cdot + x_i\big) 
- z_{K(0)}^{a,b}(\cdot)\|_{H^{1}_a(B_{2R_i}(0))} \le \eps_i,  
\end{align}
where $H^{1}_a(B_{2R_i}(0))$ is the weighted Sobolev space $\{u:
|x|^{-a}|\n u|, |x|^{-a} u\in L^2(B_{2R_i}(0))\}$.
\end{proposition}
\begin{proof}
Consider 
\begin{align*}
\phi_i(x) = u_i(x_i)^{-1} u_i\big(u_i(x_i)^{-\frac{2}{N-2-2a}}x + x_i\big), \quad
|x|<\bar r u_i(x_i)^{\frac{2}{N-2-2a}} . 
\end{align*}
We have 
\begin{align*}
-\div\Big(\big|x+&u_i(x_i)^{\frac{2}{N-2-2a}}x_i\big|^{-2a} \nabla \phi_i(x)\Big) \\
&= K_i\big(u_i(x_i)^{-\frac{2}{N-2-2a}}x + x_i\big) 
~\big|x+u_i(x_i)^{\frac{2}{N-2-2a}}x_i\big|^{-bp}~ \phi_i^{p-1}(x).    
\end{align*}
Moreover, from the definition of isolated blow-up
\begin{align}\label{eq:13}
\phi_i(0)=1,\quad 
0< \phi_i(x) \le \bar{C} |x|^{-\frac{N-2-2a}{2}} \quad\text{ for }
|x|< \bar{r}\, u_i(x_i)^{\frac{2}{N-2-2a}}.   
\end{align}
Lemma \ref{l:harnack:zero} shows that for large $i$ and for any $0<r<1$ we have 
\begin{align}
\label{eq:14}
\max_{\partial B_r}\phi_i \le C \min_{\partial B_r}\phi_i,  
\end{align}
where $C=C(N,\bar C,A_1)$. Since 
\[-\div\Big(\big|x+ u_i(x_i)^{\frac{2}{N-2-2a}}x_i\big|^{-2a} \nabla
\phi_i(x)\Big) \ge 0 \text{ and } \phi_i(0)=1\] 
we may use (\ref{eq:14}) and the minimum principle for $|x|^{-2a}$-superharmonic
functions in \cite[Thm 7.12]{potential}
 to deduce that
\begin{align}
\label{eq:3}
\phi_i(x) \le C \quad \text{ in } B_1(0).  
\end{align}
From (\ref{eq:13}), (\ref{eq:3}) and regularity results in
\cite{SF_reg} the functions $\phi_i$ are uniformly 
bounded in $C_{loc}^{0,\gamma}(\rz^N)$ and $H^1_{a,loc}(\rz^N)$ for
some $\gamma\in(0,1)$.   
Since point-concentration is ruled out by the $L^\infty$-bound, there
is some positive function $\phi 
\in C^{0,\gamma'}_{loc}(\rz^N)\cap H^1_{a,loc}(\rz^N)$ and
some $\gamma'\in(0,1)$ such that 
\begin{align*}
&\phi_i \to \phi \text{ in } C_{loc}^{0,\gamma'}(\rz^N)\cap
H^1_{a,loc}(\rz^N),\\ 
&-\div(|x|^{-2a} \nabla \phi) = \lim_{i \to \infty} K_i(x_i)
\frac{\phi^{p-1}}{|x|^{bp}}\\
&\phi(0)=1.
\end{align*}
By uniqueness of the solutions proved in \cite{ChouChu} we deduce
that $\phi=z_{K(0)}^{a,b}$.  
\end{proof}

\begin{remark}\label{r:bubble0} From the proof of Proposition
  \ref{prop:Bubble:zero}  one can easily
  check that if $x_i \to 0$ is an
  isolated blow-up point then there exists a positive constant $C$,
  depending on $\lim_{i \to 
    \infty}K_i(x_i)$ and $a$, $b$, and $N$, such that the function
  $\bar w_i$ defined in 
  (\ref{eq:62}) is strictly decreasing for $C
  u_i(x_i)^{-2/(N-2-2a)}\le r\le r_i=R_i u_i(x_i)^{-\frac
    2{N-2-2a}}$ (see Lemma \ref{l:decrease} of
  the Appendix).   
\end{remark}
\begin{lemma}\label{l:1}
Let $x_i\to 0$ be a blow-up point. Then for any $x$ such
that $|x-x_i|\geq r_i=R_i u_i(x_i)^{-\frac{2}{N-2-2a}}$ we have 
\[
|x-x_i|=|x|(1+\smallo  ).
\]
In particular, $x_i\in B_{r_i}(0)$. 
\end{lemma}
\begin{proof}
The assumption $\Big|x_i\ui(x_i)^{\frac{2}{N-2-2a}}\Big|=\smallo $ implies that
$|x_i|=r_i\,\smallo  $. Hence
\[
|x|\geq|x-x_i|-|x_i|\geq r_i-r_i\,\smallo  =r_i(1+\smallo  ).
\]
Therefore
\[
\frac{|x_i|}{|x|}\leq\frac{r_i\,\smallo  }{r_i(1+\smallo  )}=\smallo  
\]
and hence
\[
\bigg|\frac{x-x_i}{|x|}\bigg|=\bigg|\frac
x{|x|}-\frac{x_i}{|x|}\bigg|\
\mathop{\longrightarrow}\limits_{i\to+\infty}\ 1
\]
thus proving the lemma.
\end{proof}

\begin{proposition}\label{p:upperbound}
Suppose $\{K_i\}_i\subset C^1_{\rm loc}(B_2)$ satisfy (\ref{eq:blow1})
with $\Omega=B_2$ and 
\begin{align}
\label{eq:5}
|\nabla K_i(x)| \le A_2 \text{ for all } x \in B_2. 
\end{align}
Let $\ui$ satisfy $(P_i)$ with $\Omega=B_2$ and suppose that $x_i\to 0$ is
an isolated simple blow-up point such that
\begin{align}
\label{eq:10}
|x-x_i|^{\frac{N-2-2a}2} u_i(x) \le A_3 \text{ for all } x \in B_2.  
\end{align}
Then there exists $C=C(N,a,b,A_1,A_2,A_3,\bar C,\rho)>0$ such that
\begin{equation}\label{eq:upperbound}
\ui(x)\leq C\,\ui(x_i)^{-1}|x-x_i|^{2+2a-N}\quad\text{for all}\
|x-x_i|\leq1.
\end{equation}
Furthermore there exists a H\"older continuous function $B(x)$ (smooth outside $0$) satisfying
$\div(|x|^{-2a}\n B)=0$ in $B_1$, such that, after passing to a subsequence,
\[
\ui(x_i)\ui(x)\to h(x)=A|x|^{2+2a-N}+B(x)\quad\text{in}\ C^2_{\rm
  loc}(B_1\setminus\{0\})
\]
where 
\begin{align*}
A=\frac{K(0)}{(N-2-2a)|{\mathbb S}^N|} 
\int_{\R^N}\frac{\big(z_{K(0)}^{a,b}\big)^{p-1}}{|x|^{bp}}\,dx.
\end{align*}
\end{proposition}

\begin{lemma}
\label{l:yy:2.2:zero}
Under the assumption of Proposition \ref{p:upperbound} without (\ref{eq:5}) there exist
a positive $\delta_i= O\Big(R_i^{\frac{-2(1+a-b)(N-2-2a)}{N-2(1+a-b)}}\Big)$ and
$C=C(N,a,b,A_1,A_2,\bar C,\rho)>0$ such that
\begin{align}
\label{eq:20}
\ui(x)\leq C\,\ui(x_i)^{-\lambda_i}|x-x_i|^{2+2a-N+\delta_i}\quad\text{for all}\
r_i \le |x-x_i|\leq1,  
\end{align}
where $\lambda_i:= 1-2\delta_i/(N-2-2a)$.
\end{lemma}
\begin{proof}
It follows from Proposition \ref{prop:Bubble:zero} that
\begin{align}
\label{eq:1}
u_i(x) \le c u_i(x_i) R_i^{2a+2-N} \text{ for } |x-x_i|=r_i.  
\end{align}
From the definition of isolated simple blow-up in (\ref{eq:isolated-simple}) there exists $\rho>0$
such that
\begin{align}
\label{eq:2}
r^{\frac{N-2-2a}{2}} \bar u_i \text{ is strictly decreasing in } r_i<r<\rho.  
\end{align}
From (\ref{eq:1}), (\ref{eq:2}) and Lemma \ref{l:harnack:zero} it
follows that for all $r_i \le |x-x_i|<\rho$
\begin{align*}
|x-x_i|^{\frac{N-2-2a}{2}} u_i(x) &\le c  |x-x_i|^{\frac{N-2-2a}{2}} \bar u_i(|x-x_i|)
\le c r_i^{\frac{N-2-2a}{2}} \bar u_i(r_i)
\le c R_i^{\frac{2+2a-N}{2}}.
\end{align*}
Therefore for $r_i< |x-x_i|<\rho$
\begin{align}
\label{eq:21}
u_i(x)^{\frac{4}{N-2-2a}} &\le c R_i^{-2} |x-x_i|^{-2}.  
\end{align}
Consider the following degenerated elliptic operator
\begin{align*}
\mathcal{L}_i \varphi = \div(|x|^{-2a} \nabla \varphi) +   K_i(x)|x|^{-bp} u_i(x)^{p-2}\varphi.    
\end{align*}
Clearly $u_i>0$ solves $\mathcal{L}_iu_i=0$. Hence $-\mathcal{L}_i$ is nonnegative and the maximum principle 
holds for $\mathcal{L}_i$.  
Direct computations show for any $0 \le \mu \le N-2-2a$
\begin{align}
\label{eq:22}
\div(|x|^{-2a} \nabla (|x|^{-\mu})) = -\mu(N-2-2a-\mu) |x|^{-2-2a-\mu} \text{ for } x\neq 0.  
\end{align}
From (\ref{eq:21}), (\ref{eq:22}) and Lemma \ref{l:1} we infer
\begin{align*}
\mathcal{L}_i(|x|^{-\mu})\le \big(-\mu(N-2-2a-\mu)+c R_i^{\frac{-2(1+a-b)(N-2-2a)}{N-2(1+a-b)}}\big) |x|^{-2-2a-\mu} .
\end{align*}
We can choose $\delta_i= O( R_i^{\frac{-2(1+a-b)(N-2-2a)}{N-2(1+a-b)}})$ such that
\begin{align}
\label{eq:23}
\max(\mathcal{L}_i(|x|^{-\delta_i}),\mathcal{L}_i(|x|^{2a+2-N+\delta_i})) \le 0. 
\end{align}
Set $M_i:= 2\max_{\partial B_\rho(x_i)}u_i$, $\lambda_i=1-2\delta_i/(N-2-2a)$, and
\begin{align}
\label{eq:25}
\phi_i(x):= M_i \rho^{\delta_i}|x|^{-\delta_i}+ A u_i(x_i)^{-\lambda_i} |x|^{2+2a-N+\delta_i} 
\text{ for } r_i\le |x-x_i|\le \rho,
\end{align}
where $A$ will be chosen later. We will apply the maximum principle to compare $\phi_i$ and $u_i$.
By the choice of $M_i$ and  Lemma \ref{l:1} we infer for $i$ sufficiently large
\begin{align*}
\phi_i(x) \ge \frac{M_i}2 \ge u_i(x) \text{ for } |x-x_i|=\rho.  
\end{align*}
On the inner boundary $|x-x_i|=r_i$ we have by (\ref{eq:1}) and for $A$ large enough
\begin{align*}
\phi_i(x) &\ge A (1+o(1))u_i(x_i)^{-\lambda_i}r_i^{2+2a-N+\delta_i}= A
(1+o(1))R_i^{2+2a-N+\delta_i} 
u_i(x_i)^{2-\frac{2\delta_i}{N-2-2a}-\lambda_i}\\
&\ge  A (1+o(1)) R_i^{2+2a-N} u_i(x_i)  \ge u_i(x).  
\end{align*}
Now we obtain from the maximum principle in the annulus $r_i\le |x-x_i| \le \rho$ that
\begin{align}
\label{eq:24}
u_i(x) \le \phi_i(x) \text{ for all } r_i \le |x-x_i|\le \rho.  
\end{align}
It follows from (\ref{eq:2}), (\ref{eq:24}) and Lemma \ref{l:harnack:zero} that for any $r_i\le
\theta \le \rho$ we have
\begin{align*}
\rho^{\frac{N-2-2a}{2}} M_i &\le c  \rho^{\frac{N-2-2a}{2}} \bar u_i(\rho)
\le c \theta^{\frac{N-2-2a}{2}} \bar u_i(\theta) \\
&\le c \theta^{\frac{N-2-2a}{2}} \big(M_i\rho^{\delta_i}\theta^{-\delta_i}+A u_i(x_i)^{-\lambda_i}\theta^{2+2a-N+\delta_i}\big). 
\end{align*}
Choose $\theta=\theta(\rho,c)$ such that 
\begin{align*}
 c \theta^{\frac{N-2-2a}{2}} \rho^{\delta_i}\theta^{-\delta_i}<\frac{1}{2}\rho^{\frac{N-2-2a}{2}}.
\end{align*}
Then we have 
\begin{align*}
 M_i \le c u_i(x_i)^{-\lambda_i}, 
\end{align*}
which, in view of (\ref{eq:24}) and the definition of $\phi_i$ in (\ref{eq:25}), proves
(\ref{eq:20}) for $r_i\le |x-x_i|\le \rho$. The Harnack inequality in Lemma
\ref{l:harnack:zero} allows to extend (\ref{eq:20}) for $r_i\le |x-x_i|\le 1$.
\end{proof}

\begin{proof}[Proof of Proposition \ref{p:upperbound}]
The inequality of Proposition \ref{p:upperbound} for $|x-x_i|\leq r_i$
follows immediately for Proposition \ref{prop:Bubble:zero}. Let $e\in\R^N$,
$|e|=1$ and consider the function
\[
v_i(x)=\ui(x_i+e)^{-1}\ui(x).
\]
Clearly $v_i$ satisfies the equation
\begin{equation}\label{eq:vi}
-\div (|x|^{-2a}\n
 v_i)=\ui(x_i+e)^{p-2}K_i(x)\frac{{v_i}^{p-1}}{|x|^{bp}}\quad\text{in}\
 B_{4/3}.
\end{equation}
Applying the Harnack inequality of Lemma \ref{l:harnack:zero} on
$v_i$, we obtain that $v_i$ is bounded on any compact set not
containing $0$. By standard elliptic theories, it follows that, up to
a subsequence, $\{v_i\}_i$ converges in $C^2_{\rm
  loc}(B_2\setminus\{0\})$ to some positive function $v\in
C^2(B_2\setminus\{0\})$. Since $\ui(x_i+e)\to 0$ due to Lemma
\ref{l:yy:2.2:zero},
 we can pass to the limit in (\ref{eq:vi}) thus obtaining 
\[
-\div(|x|^{-2a}\n v)=0\quad\text{in}\ B_2\setminus\{0\}.
\]
We claim that $v$ has a singularity at $0$. Indeed, from Lemma
\ref{l:harnack:zero}
 and standard elliptic theories, for any $0<r<2$ we have that
\[
\lim_{i\to\infty}\ui(x_i+e)^{-1}r^{\frac{N-2-2a}2}\bar\ui(r)=r^{\frac{N-2-2a}2}\bar
v(r)
\]
where $\bar v(r)=\frac1{|\partial B_r|}\int_{\partial B_r}v$. Since
the blow-up is simple isolated, $r^{\frac{N-2-2a}2}\bar v(r)$ is
non-increasing for $0<r<\rho$ and this would be impossible in the case
in which $v$ is regular at $0$. It follows that $v$ is singular at $0$
and hence from the B\^ocher-type Theorem proved in the Appendix (see
Theorem \ref{t:bocher})
\[
v(x)=a_1|x|^{2+2a-N}+b_1(x)
\]
where $a_1>0$ is some positive constant and $b_1(x)$ is some H\"older 
continuous function in $B_2$ such that $-\div(|x|^{-2a}\n b_1)=0$. 

Let us first establish the inequality in Proposition
\ref{p:upperbound} for $|x-x_i|=1$. Namely we prove that
\begin{equation}\label{eq:mod1}
\ui(x_i+e)\leq c\,\ui(x_i)^{-1}.
\end{equation}
By contradiction, suppose that (\ref{eq:mod1}) fails. Then along a
subsequence, we have 
\begin{equation}\label{eq:mod11}
\lim_{i\to\infty}\ui(x_i+e)\ui(x_i)=\infty.
\end{equation}
Multiplying $(P_i)$ by $\ui(x_i+e)^{-1}$ and integrating on $B_1$, we
get
\begin{equation}\label{eq:mod12}
-\int_{\partial B_1}|x|^{-2a}\frac{\partial v_i}{\partial
 \nu}= \int_{B_1}K_i(x)\frac{\ui^{p-1}}{|x|^{bp}}\ui(x_i+e)^{-1}\,dx.
\end{equation}
From the properties of $b_1$ and the convergence of $v_i$ to $v$, we
know that 
\begin{align}\label{eq:mod13}
\lim_{i\to\infty}\int_{\partial B_1}|x|^{-2a}\frac{\partial
  v_i}{\partial \nu}&=\int_{\partial
  B_1}|x|^{-2a}\frac{\partial}{\partial
  \nu}\big(a_1|x|^{2+2a-N}+b_1(x)\big)\notag\\
&=-a_1(N-2-2a)|{\mathbb S}^N|<0.
\end{align} 
From Proposition \ref{prop:Bubble:zero} there holds
\begin{equation}\label{eq:mod14}
\int_{|x-x_i|\leq r_i}\frac{K_i(x)\ui^{p-1}}{|x|^{bp}}\,dx\leq
C\,\ui(x_i)^{-1}
\end{equation}
while from Lemma \ref{l:yy:2.2:zero} and Lemma \ref{l:1} we have that 
\begin{align}\label{eq:mod15}
\int_{r_i\leq |x-x_i|\leq 1}\frac{K_i(x)\ui^{p-1}}{|x|^{bp}}\,dx&\leq
c\int_{r_i\leq |x-x_i|\leq
  1}\ui(x_i)^{-\l_i(p-1)}\frac{|x-x_i|^{(2+2a-N+\d_i)(p-1)}}{|x|^{bp}}\notag\\
&\leq c\,\ui(x_i)^{-\l_i(p-1)}r_i^{(2+2a-N+\d_i)(p-1)-bp+N}\notag\\
&=c\,\ui(x_i)^{-1}R_i^{(2+2a-N+\d_i)(p-1)-bp+N}=o(1)\ui(x_i)^{-1}.
\end{align}
Finally, (\ref{eq:mod11}), (\ref{eq:mod13}), (\ref{eq:mod14}), and
(\ref{eq:mod15}) lead to a contradiction. Since we have established
(\ref{eq:mod1}), the inequality in Proposition \ref{p:upperbound} has
been established for $\rho\leq|x-x_i|\leq1$ (due to Lemma  \ref{l:harnack:zero}). It
remains to treat the case $r_i\leq|x-x_i|\leq\rho$. To 
this aim we scale the problem to reduce it to the case $|x-x_i|=1$. By
contradiction, suppose that there exists a subsequence $\tilde x_i$
satisfying $r_i\leq |\tilde x_i-x_i|\leq \rho$ and 
\begin{equation}\label{eq:mod16}
\lim_{i\to+\infty}\ui(\tilde x_i)\ui(x_i)|\tilde
x_i-x_i|^{N-2-2a}=+\infty.
\end{equation}
Set $\tilde r_i=|\tilde x_i-x_i|$ and $\tilde \ui(x)=\tilde
r_i^{\frac{N-2-2a}2}\ui(\tilde r_i x)$. We have that $\tilde\ui$ satisfies the equation
\[
-\div(|x|^{-2a}\n \tilde \ui(x))=K_i(\tilde r_i
 x)\frac{\tilde\ui(x)^{p-1}}{|x|^{bp}}.
\]
Since $|x_i|=r_i\,\smallo$ and $\tilde r_i\geq r_i$ we have that
$x_i/\tilde r_i\to 0$. We have that $x_i/\tilde r_i$ is an isolated simple blow-up point for
$\{\tilde \ui\}_i$. From (\ref{eq:mod1}), we have that 
\[
\tilde\ui\bigg(\frac {x_i}{\tilde r_i}+\frac {\tilde x_i-x_i}{\tilde
  r_i}\bigg)\leq c\tilde \ui\bigg(\frac {x_i}{\tilde r_i}\bigg)^{-1}
\]
which gives
\[
\tilde r_i^{N-2-2a}\ui(\tilde x_i)\ui(x_i)\leq c.
\]
The above estimate and (\ref{eq:mod16}) give rise to a
contradiction. The inequality in Proposition \ref{p:upperbound} is
thereby established. \par
We compute $A$ by multiplying $(P_i)$ by $\ui(x_i)$ and
integrating over $B_1$. From the divergence theorem,
\begin{equation}\label{eq:mod17}
-\int_{\partial
 B_1}|x|^{-2a}\frac{\partial}{\partial\nu}\big(\ui(x_i)u_i\big)= \ui(x_i)\int_{B_1}K_i(x)\frac{\ui^{p-1}}{|x|^{bp}}\,dx.
\end{equation}
Let $w_i(x)=\ui(x_i)u(x)$. We have that $w_i$ satisfies
\[
-\div(|x|^{-2a}\n w_i)=\ui(x_i)^{2-p}K_i(x)\frac{w_i^{p-1}}{|x|^{bp}}.
\]
Moreover the inequality (\ref{eq:upperbound}) implies that $w_i$ is
bounded on any compact set not containing $0$. Hence $w_i\to w$ in
$C^2_{\rm loc}(\R^N\setminus\{0\})$ where $w$ satisfies
\[
-\div(|x|^{-2a}\n w)=0\quad\text{in}\ \R^N\setminus\{0\}.
\]
From the B\^ocher-type theorem proved in the Appendix (Theorem
\ref{t:bocher}), we find that $w(x)=A|x|^{2+2a-N}+B(x)$ 
where $B(x)$ is H\"older continuous in $\R^N$ and satisfies 
$ -\div(|x|^{-2a}\n B)=0\quad\text{in}\ \R^N$. Hence 
\begin{align}\label{eq:mod18}
\lim_{i\to\infty}\int_{\partial B_1}|x|^{-2a}\frac{\partial}{\partial
  \nu}\big(\ui(x_i)\ui\big)&=\int_{\partial B_1}|x|^{-2a}\frac{\partial}{\partial
  \nu}\big(A|x|^{2+2a-N}+B(x)\big)\notag\\
&=A(2+2a-N)|{\mathbb S}^N|.
\end{align}
On the other hand from (\ref{eq:mod15}) and Proposition \ref{prop:Bubble:zero}
\begin{align}\label{eq:mod19}
\ui(x_i)\int_{B_1}K_i(x)\frac{\ui^{p-1}}{|x|^{bp}}\,dx &=  \ui(x_i)\int_{|x-x_i|\leq
  r_i}K_i(x)\frac{\ui^{p-1}}{|x|^{bp}}\,dx + \smallo \notag\\
&=  K_i(0)\int_{|y|\leq R_i}\frac{\big(z_{K(0)}^{a,b}\big)^{p-1}}{\big|y+\ui(x_i)^{\frac2{N-2-2a}}x_i\big|^{bp}}\,dy 
+ \smallo \notag\\
&= K(0)\int_{\R^N}\frac{\big(z_{K(0)}^{a,b}\big)^{p-1}}{|y|^{bp}}\,dy+o(1).
\end{align}
By (\ref{eq:mod17}), (\ref{eq:mod18}), and (\ref{eq:mod19}) the value of $A$ is computed and
Proposition \ref{p:upperbound} is thereby established.
\end{proof}

Using Proposition \ref{prop:Bubble:zero} and the upper bound in Proposition \ref{p:upperbound} it is
easy to see that the following estimates hold. 
\begin{lemma}
\label{l:yy:2.4:zero}
Under the assumptions of Proposition \ref{p:upperbound} we have for $s=s_1+s_2$
\begin{align*}
&\int_{|x-x_i|\le r_i} |x-x_i|^{s_1}|x|^{s_2} |x|^{-bp} u_i(x)^p \\
&\quad = 
\begin{cases}
u_i(x_i)^{\frac{-2s}{N-2-2a}} \left(o(1)+\int_{\rz^N}|x|^{s-bp} z_{1,K_i(x_i)}^p \right) 
&\text{ if } -N+bp<s<N-bp,\\
O(u_i(x_i)^{-p} \log u_i(x_i)) &\text{ if } s=N-bp,\\
o\big(u_i(x_i)^{-p}\big)  &\text{ if } s>N-bp.
\end{cases}
\\
&\int_{r_i\le |x-x_i|\le 1} |x-x_i|^{s_1}|x|^{s_2}  |x|^{-bp} u_i(x)^p \\
&\quad \le 
\begin{cases}
o\big(u_i(x_i)^{\frac{-2s}{N-2-2a}}\big) &\text{ if } -N+bp<s<N-bp,\\
O(u_i(x_i)^{-p} \log u_i(x_i)) &\text{ if } s=N-bp,\\
O\big(u_i(x_i)^{-p}\big)  &\text{ if } s>N-bp.
\end{cases}
\end{align*}

\end{lemma}

\begin{proposition}\label{p:simplezero}
Let $a\in\big[\frac{N-4}2,\frac{N-2}2\big[$. Suppose
 that $\{K_i\}_i$ satisfy (\ref{eq:blow1}) with
$\Omega=B_2\subset\R^N$ for some positive constant $A_1$,
$\n K_i(0)=0$, $\{K_i\}_i$ converge to $K$ in $C^2(B_2)$, 
$\{u_i\}_i$ satisfy $(P_i)$ with $\Omega=B_2(0)$ and $x_i\to 0$ is an isolated blow-up point
with (\ref{eq:10}) for some positive constant $A_3$. Then it has to be
an isolated simple blow-up point.
\end{proposition}
\begin{proof}
  From Remark \ref{r:bubble0} there exists a constant $c$ such that
$r^{\frac{N-2-2a}2}\bar u_i(r)$ is decreasing in $c u_i(x_i)^{-2/(N-2-2a)}\le r\le r_i$. 
Arguing by contradiction, let us suppose that the blow-up is not
simple. Hence for any $i$ there exists $\mu_i\geq r_i$, $\mu_i\to0$,
such that $\mu_i$ is the first point after $r_i$ in which the function
$r^{\frac{N-2-2a}2}\bar u_i(r)$ becomes increasing. In particular
$\mu_i$ is a critical point of such a function. Set 
\[
\xi_i(x)=\mu_i^{\frac{N-2-2a}2}u_i(\mu_i x),\quad\text{for}\
|\mu_i x-x_i|\leq1.
\]
Clearly $\xi_i$ satisfies
\[
-\div(|x|^{-2a}\n \xi_i)=  K_i(\mu_i
x)\frac{\xi_i^{p-1}}{|x|^{bp}},
\quad\text{for}\
|\mu_i x-x_i|\leq1.
\]
Note that $\mu_i^{-1}\leq
R_i^{-1}u_i(x_i)^{\frac2{N-2-2a}}\leq u_i(x_i)^{\frac2{N-2-2a}}$ and
hence 
$$
\mu_i^{-1}|x_i|\leq u_i(x_i)^{\frac2{N-2-2a}}|x_i|\to~0
$$ 
in view of (\ref{eq:16}). Moreover (\ref{eq:10}) implies that 
\[
|x-\mu_i^{-1}x_i|^{\frac{N-2-2a}2}\xi_i(x)\leq {\rm const}\quad\text{for}\
|x-\mu_i^{-1}x_i|\leq1/\mu_i.
\]
It is also easy to verify that 
\[
\lim_{i\to\infty}\xi_i(\mu_i^{-1}x_i)=\lim_{i\to\infty}\mu_i^{\frac{N-2-2a}2}u_i(x_i)
=\infty.
\]
On the other hand 
\begin{align*}
\media_{\partial B_r(\mu_i^{-1}x_i)}\xi_i=\mu_i^{\frac{N-2-2a}2}\media_{\partial
  B_{r\mu_i}(x_i)}u_i=\mu_i^{\frac{N-2-2a}2}\bar u_i(\mu_i r).
\end{align*}
Hence 
\[
r^{\frac{N-2-2a}2}\bar\xi_i(r)=\bar w_i(\mu_i r)
\]
and the function
$r^{\frac{N-2-2a}2}\bar\xi_i(r)$ is decreasing in
$c\xi_i(\mu_i^{-1}x_i)^{-\frac2{N-2-2a}}<r<1$ so that $0$ is an isolated simple
blow-up point for $\{\xi_i\}$. From Proposition \ref{p:upperbound} we
have that 
\[
\xi_i(\mu_i^{-1}x_i)\xi_i(x)\to h(x)=A|x|^{2+2a-N}+B(x)\quad\text{in }C^2_{\rm
  loc}(\R^N\setminus\{0\})
\]
where  $B(x)$ is H\"older continuous in $\R^N$ and satisfies 
$ -\div(|x|^{-2a}\n B)=0\quad\text{in}\ \R^N$. Since $h\geq0$, the Harnack
inequality implies that $B$ is bounded and from the Liouville Theorem
(see \cite{potential}) we find that $B$ must be constant. Since 
\[
\frac d{dr}\{h(r)r^{\frac{N-2-2a}2}\}|_{r=1}=0
\]
we have that $A=B>0$. 
From the Taylor expansion, (\ref{eq:16}) and the assumption
$\n K_i(0)=0$ we find
\begin{align}\label{eq:54}
|\n K_i(\mu_i^{-1}x_i)|\leq {\rm const} |\mu_i^{-1}x_i|=o\Big(\xi_i(\mu_i^{-1}x_i)^{-\frac2{N-2-2a}}\Big).
\end{align}
Using Lemma \ref{l:yy:2.4:zero}, \eqref{eq:54}, and the assumption
on $a$, we have 
\begin{align*}
&\int_{B_{\sigma}(\mu_i^{-1}x_i)}\mu_i \n K_i(\mu_i
x)\cdot
x\frac{\xi_i^p}{|x|^{bp}}=\int_{B_{\sigma}(\mu_i^{-1}x_i)}\mu_i\Big[\n
K_i(\mu_i^{-1}x_i)+O(\mu_i x-\mu_i^{-1}x_i)\Big]\cdot
x\frac{\xi_i^p}{|x|^{bp}}\\
&=\int_{B_{\sigma}(\mu_i^{-1}x_i)}\mu_i\Big[\n
K_i(\mu_i^{-1}x_i)+O(|x|+|x-\mu_i^{-1}x_i|)\Big]\cdot
x\frac{\xi_i^p}{|x|^{bp}}
\\
&\quad=
\left.
\begin{cases}
\mu_i O\Big(\xi_i(\mu_i^{-1}x_i)^{-\frac4{N-2-2a}}\Big)&\text{if } p>\frac4{N-2-2a}\\
\mu_i O\Big(\xi_i(\mu_i^{-1}x_i)^{-p}\log u_i(x_i)\Big)&\text{if }p=\frac4{N-2-2a}\\
\mu_i O\Big(\xi_i(\mu_i^{-1}x_i)^{-p}\Big)&\text{if }p<\frac4{N-2-2a}
\end{cases}
\right\}
=o(\xi_i(\mu_i^{-1}x_i)^{-2}).
\end{align*}
Hence, from Corollary \ref{c:poho1} and (\ref{eq:upperbound}), 
we have that for any $0<\sigma<1$
\begin{align*}
&\int_{\partial B_{\sigma}(0)}B(\sigma,x,\xi_i,\n
\xi_i)\\
&\quad =\frac{1}{p}\int_{B_{\sigma}(0)}\mu_i \n K_i(\mu_i x)\cdot x\frac{\xi_i^p}{|x|^{bp}}
- \frac{\sigma}{p}\int_{\partial B_{\sigma}(0)}K_i(\mu_i x)\frac{\xi_i^p}{|x|^{bp}}\\
&\quad =\frac{1}{p}\int_{B_{\sigma}(\mu_i^{-1}x_i)}\mu_i \n K_i(\mu_i x)\cdot x\frac{\xi_i^p}{|x|^{bp}}
+O(\xi_i(\mu_i^{-1}x_i)^{-p})\\
&\quad =o(\xi_i(\mu_i^{-1}x_i)^{-2}).
\end{align*}
Multiplying by $\xi_i(\mu_i^{-1}x_i)^{2}$ and letting $i\to\infty$ we find that
\[
\int_{\partial B_{\sigma}}B(\sigma,x,h,\n h)=0.
\]
On the other hand Proposition \ref{p:boundterm} implies that for small
$\sigma$ the above integral is strictly negative, thus giving rise to a
contradiction. The proof is now complete.
\end{proof}

\section{A-priori estimates}
To prove the a-priori estimates we first locate the possible blow-up points as in 
\cite{SchoenZhang}. To this end we use the Kelvin transform defined in (\ref{eq:kelvin}).
We recall that if $u$ solves (\ref{eq:49}) then $\tilde u= |x|^{-(N-2-2a)}u(x/|x|^2)$ solves (\ref{eq:49})
with $K$ replaced by $\tilde K(x)=K(x/|x|^2)$. Since weak solutions to (\ref{eq:49}) are H\"older continuous
(see \cite{SF_reg}) we infer that 
\begin{align}
\label{eq:63}
\lim_{|x|\to \infty} |x|^{N-2-2a}u(x) \text{ exists.}  
\end{align}
Let us define $\omega_a(x):=(1+|x|^{N-2-2a})^{-1}$. 
\begin{lemma}
\label{l:standardblowup}
Suppose $a\ge 0$, $2<p<2^*$, and $K\in C^2 (\rz^N)$ satisfies (\ref{eq:60}) and for some positive
constants $A_1$, $A_2$ condition (\ref{eq:51}) and 
\begin{align}
\label{eq:61}
\|\n K\|_{L^{\infty}(B_2(0))}+ \|\n \tilde K\|_{L^{\infty}(B_2(0))}\leq A_2.  
\end{align}
Then for any $\e\in(0,1)$, $R>1$, there exists
$C_0 =C_0(\e,R,N,a,b,A_1,A_2)>0$, 
such that if $u$ is a solution of
(\ref{eq:49}) and ${\mathcal K}=\{q_1,\dots ,q_k\}\subset \R^N\cup \{\infty\}$ with 
\begin{align}
\label{eq:31} \raisetag{45pt}
\begin{cases}
\max\limits _{x\in \R^N\cup\{\infty\}}\frac{u(x)}{\omega_a(x)}\dist(x,{\mathcal K})^{\frac{N-2-2a}{2}}>C_0,\\[5pt]
u(q_i)|q_i|^{\frac{2}{N-2-2a}}<\eps,\text{ and for all } 1\le i \le k\\[5pt]
\max\limits _{x\in \R^N\cup\{\infty\}}\frac{u(x)}{\omega_a(x)}\dist(x,\{q_1,\dots,q_{i-1}\})^{\frac{N-2-2a}{2}}
\le  \frac{u(q_i)}{\omega_a(q_i)}\dist(q_i,\{q_1,\dots,q_{i-1}\})^{\frac{N-2-2a}{2}},     
\end{cases}
\end{align}
then there exists $q^* \not\in {\mathcal K}$ such that $q^*$
is a maximum point of $({u}/{\omega_a})\dist(\cdot,{\mathcal
  K})^{\frac{N-2-2a}{2}}$~and
\begin{itemize}
\item[(A)] if $|q^*|\le 1$
\begin{align}\label{eq:66}
\bigg\|\dfrac{u(u(q^*)^{-\frac{2}{N-2-2a}}x+q^*)}{u(q^*)}-z^{a,b}_{K(q^*)}\bigg\|_{\substack{C^{0,\gamma}
(B_{2R}(0))\cap\\ H^1_a(B_{2R}(0))}}+ 
|q^*|u(q^*)^{\frac{2}{N-2-2a}}<\e
\end{align}
\item[(B)] if $|q^*|>1$
\begin{align}\label{eq:67}
\bigg\|\dfrac{ \tilde u(\tilde u(\tilde
 q^*)^{-\frac{2}{N-2-2a}}x+\tilde q^*)}{\tilde u(\tilde
 q^*)}-z^{a,b}_{K( q^*)}\bigg\|_{\substack{C^{0,\gamma} 
(B_{2R}(0))\cap\\ H^1_a(B_{2R}(0))}}+
|\tilde q^*|\tilde u(\tilde q^*)^{\frac{2}{N-2-2a}}<\e
\end{align}
\end{itemize}
where $\tilde q^* = {\text Inv}(q^*):=q^*/|q^*|^2$, $\tilde u$ is the Kelvin transform of $u$, 
$\dist(\cdot,\cdot)$ is the
distance on $\rz^N\cup \{\infty\}$ induced by the standard metric 
on the sphere through the stereo-graphic projection, and $\dist(\cdot,\emptyset)\equiv 1$.
\end{lemma}
\begin{proof}
Fix $\eps>0$ and $R>1$. Let $C_0$ and $C_1$ be  positive constants
depending on $\e,R,a,b,N, A_1,A_2$ which shall be appropriately chosen
in the sequel.\\
Let $q^*\in \rz^N \cup \{\infty\}$ be the maximum point of 
$u/\omega_a \dist(x,{\mathcal K})^{\frac{N-2-2a}{2}}$. By (\ref{eq:63}) this maximum is 
achieved. From the first in (\ref{eq:31}) we have that
$u(q^*)/\omega_a(q^*) \dist(q^*,{\mathcal K})^{\frac{N-2-2a}{2}}>C_0$.
First we treat the case  $|q^*|\leq 1$. We claim that there exists a
constant $C_1$, depending only on  $\e,R,a,b,N, A_1,A_2$, such that 
$|q^*|^{\frac{N-2-2a}2}u(q^*)< C_1$. 
If not, there exist
solutions $u_i$ of (\ref{eq:49}) and finite sets ${\mathcal K}_i=\{q_1^i,\dots, q_{k_i}^i\}$ satisfying
(\ref{eq:31}) above, such that  
for the maximum points $q^*_i$ of $u_i/\omega_a \dist(\cdot,{\mathcal
  K}_i)^{\frac{N-2-2a}{2}}$ there holds
\[|q^*_i|\leq1 \text{ and } |q^*_i|^{\frac{N-2-2a}2}u_i(q_i^*)\to
\infty.\]
Consider the functions $v_i$, defined by
\begin{align*}
v_i(x) := u_i(q_i^*)^{-1}u_i\big(|q_i^*|^{1+\frac{(N-2-2a)(2-p)}{4}}u(q_i^*)^{\frac{2-p}{2}}x+q_i^*\big),  
\end{align*}
which satisfy
\begin{align*}
&-\div\big(\big||q_i^*|^{\frac{(N-2-2a)(2-p)}{4}}u_i(q_i^*)^{\frac{2-p}2}x+|q_i^*|^{-1}q_i^*\big|^{-2a}
\nabla v_i\big)\\
&\quad = K\big(|q_i^*|^{1+\frac{(N-2-2a)(2-p)}{4}}u_i(q_i^*)^{\frac{2-p}2}x+q_i^*\big)\frac{v_i^{p-1}}{\big||q_i^*|^{\frac{(N-2-2a)(2-p)}{4}}u_i(q_i^*)^{\frac{2-p}2}x+|q_i^*|^{-1}q_i^*\big|^{bp}}.
\end{align*}
Let $p_i=q_{j_i}^i\in {\mathcal K}_i$ be such
that $\dist(q_i^*,{\mathcal K}_i)=\dist(q_i^*,p_i)$ and set
 $\hat {\mathcal K}_i=\{q_1^i,\dots, q^i_{j_i-1}\}$.
From (\ref{eq:31}) we infer
\begin{align*}
&\dist(p_i,\hat {\mathcal K}_i) \le \dist(p_i,q_i^*)+\dist(q_i^*,\hat {\mathcal K}_i) \le 2 \dist(q_i^*,\hat {\mathcal K}_i),\\
&u_i(p_i) |p_i|^{\frac{2}{N-2-2a}}<\eps, \quad u_i(q_i^*) 
\le u_i(p_i) \bigg(\frac{\dist(p_i,\hat {\mathcal K}_i)}{\dist(q_i^*,\hat {\mathcal K}_i)}\bigg)^{\frac{N-2-2a}{2}} 
\frac{\omega_a(q_i^*)}{\omega_a(p_i)}, 
\end{align*}
and finally that if $|p_i|\le 2$
\begin{align*}
\eps\bigg(\frac{|q_i^*|}{|p_i|}\bigg)^{\frac{2}{N-2-2a}}&>u_i(p_i) |q_i^*|^{\frac{2}{N-2-2a}}\ge
u_i(q_i^*) |q_i^*|^{\frac{2}{N-2-2a}} 
\bigg(\frac{\dist(q_i^*,\hat {\mathcal K}_i)}{\dist(p_i,\hat {\mathcal K}_i)}\bigg)^{\frac{N-2-2a}{2}} 
\frac{\omega_a(p_i)}{\omega_a(q_i^*)}\\
&\geq {\rm const}\, u_i(q_i^*) |q_i^*|^{\frac{2}{N-2-2a}} \to \infty.   
\end{align*}
Consequently there exists a positive constant $c$ such that $|q_i^*|^{-1}
\dist(q_i^*,{\mathcal K}_i)>c$, which is trivial in the case $|p_i|>2$ and follows from
the above estimate if $|p_i|\le 2$.  
Thus
\begin{align*}
|q_i^*|^{-1-\frac{(N-2-2a)(2-p)}{4}}u_i(q^*_i)^{-\frac{2-p}{2}} \dist(q_i^*,{\mathcal K}_i) &\ge
\big(u_i(q_i^*) |q_i^*|^{\frac{N-2-2a}2}\big)^{\frac{p-2}{2}} |q_i^*|^{-1}  \dist(q_i^*,{\mathcal K}_i)\\
&\ge {\rm  c} \big(u_i(q_i^*) |q_i^*|^{\frac{N-2-2a}2}\big)^{\frac{p-2}{2}} \to \infty.
\end{align*}
For $|x|\leq \frac {c}4|q_i^*|^{-\frac{(N-2-2a)(2-p)}{4}}u_i(q^*_i)^{-\frac{2-p}{2}}$ we have that
\begin{align*}
v_i(x)& =
u_i(q^*_i)^{-1}u_i\big(|q_i^*|^{1+\frac{(N-2-2a)(2-p)}{4}}u(q^*_i)^{\frac{2-p}{2}}x+q_i^*\big)\\
&\leq 
u_i(q^*_i)^{-1}\omega_a\big(|q_i^*|^{1+\frac{(N-2-2a)(2-p)}{4}}
u(q^*_i)^{\frac{2-p}{2}}x+q_i^*\big)\frac{u_i(q_i^*)}{\omega_a (q_i^*)}
\\
&\leq c\sup_{|x|\leq \frac c4}\omega_a(x+q_i^*)\omega_a(q_i^*)^{-1}\leq {\rm const}.
\end{align*}
Up to a subsequence, we have that $q_i^*\to\bar q_1$ and $v_i$ converges in $C^2_{\rm loc}(\R^N)$ to a
solution~of
\begin{align*}
-\laplace w= K(\bar q_1) w^{p-1} \text{ in }\rz^N, \quad w(0)=1.   
\end{align*}
This is impossible since the above
equation has no solution for $p<2^*$.  The claim is thereby proved.
The function $v_1$, defined by
\begin{align*}
v_1(x) := u(q^*)^{-1}u\big(u(q^*)^{-\frac{2}{N-2-2a}}x+q^*\big),  
\end{align*}
satisfies
\begin{align*}
&-\div\big(|x+u(q^*)^{\frac{2}{N-2-2a}}q^*|^{-2a} \nabla v_1\big)
= K\big(u(q^*)^{-\frac{2}{N-2-2a}}x+q^*\big)\frac{v_1^{p-1}}{|x+u(q^*)^{\frac{2}{N-2-2a}}q^*|^{bp}},\\
&v_1(0)=1.
\end{align*}
For $|x|\leq C_0^{-\frac{1}{N-2-2a}}u(q^*)^{\frac{2}{N-2-2a}} \dist(q^*,{\mathcal K})$ 
we obtain
\begin{align*}
\dist(u(q^*)^{-\frac{2}{N-2-2a}}x+q^*,{\mathcal K}) &\ge \dist(q^*,{\mathcal K}) -cC_0^{-\frac{1}{N-2-2a}} \dist(q^*,{\mathcal K})\\ 
&\ge \dist(q^*,{\mathcal K}) \big(1-cC_0^{-\frac{1}{N-2-2a}}\big)
\end{align*}
and
\begin{align*}
v_1(x)& =u(q^*)^{-1}u\big(u(q^*)^{-\frac{2}{N-2-2a}}x+q^*\big)
\\
&\leq 
u(q^*)^{-1}\omega_a\big(u(q^*)^{-\frac{2}{N-2-2a}}x+q^*\big)\frac{u(q^*)}{\omega_a (q^*)}
\left(\frac{\dist(q^*,{\mathcal K})}
{\dist(u(q^*)^{-\frac{2}{N-2-2a}}x+q^*,{\mathcal K})}\right)^{\frac{N-2-2a}{2}}
\\
&\leq \omega_a(q^*)^{-1} \big(1-cC_0^{-\frac{1}{N-2-2a}}\big)^{-\frac{N-2-2a}{2}}.
\end{align*}
Notice that $|q^*|< {\rm const\,}C_1^{\frac2 {N-2-2a}}C_0^{-\frac2
  {N-2-2a}}$ and 
$$
C_0^{-\frac{1}{N-2-2a}}u(q^*)^{\frac{2}{N-2-2a}} \dist(q^*,{\mathcal
  K})>\left(\frac{1}{4}C_0\right)^{\frac{1}{N-2-2a}}.
$$ Hence for any $\delta>0$ we may choose $C_0$,
depending on 
$a,b,N,\eps,R,A_1,A_2,C_1$, such that  
$$
\omega_a(q^*)^{-1} \big(1-C_0^{-\frac{1}{N-2-2a}}\big)^{-\frac{N-2-2a}{2}} \leq 1+\delta 
$$ 
and $v_1$ is $\eps/4$-close in $C^{0,\gamma}(B_{2R}(0))$ to a
solution of
\begin{align*}
-&\div(|x+u(q^*)^{\frac{2}{N-2-2a}}q^*|^{-2a} \nabla w)= K(q^*)
\frac{w^{p-1}}{|x+u(q^*)^{\frac{2}{N-2-2a}}q^*|^{bp}} \text{ in
  }\rz^N,\\ 
& w(0)=1,\quad 0\le w (x)\le 1+\delta.   
\end{align*}
If we choose $\delta$ small enough,
depending on $\eps$ and $R$, then it is easy to see that any solution
of the above equation is $\eps/4$-close in
$C^{0,\gamma}(B_{2R}(0))\cap H^{1}_a(B_{2R}(0))$ to
$z_{K(q^*)}^{a,b}$ and $u(q^*)^{\frac{2}{N-2-2a}}|q^*|\le
\eps/2$. This gives estimate (\ref{eq:66}). Case (B) can be reduced to
case (A) using the Kelvin transform.    
\end{proof}
\begin{proposition}\label{p:standardblowup}
Under the assumptions and notations of Lemma~\ref{l:standardblowup} there exists
for any $0<\e<1$ and $R>1$ a constant $C_0 =C_0 (\e,R,N,a,b,A_1,A_2)>0$ 
such that if $u$ is a solution of
(\ref{eq:49}) with 
\[
\max _{x\in \R^N\cup\{\infty\}}\frac{u(x)}{\omega_a(x)}>C_0 
\]
then there exist $1\leq k=k(u)<\infty$ and a set ${\mathcal
  S}(u)=\{q_1,q_2,\dots,q_k\}\subset\R^N\cup\{\infty\}$  such that
 for each $1\leq j\leq k$ we have
\begin{itemize}
\item[(A)] if $|q_j|\le 1$
\begin{align}\label{eq:64}
\bigg\|\dfrac{u(u(q_j)^{-\frac{2}{N-2-2a}}x+q_j)}{u(q_j)}-z^{a,b}_{K(q_j)}\bigg\|_{\substack{C^{0,\gamma}
(B_{2R}(0))\cap\\ H^1_a(B_{2R}(0))}}+ 
|q_j|u(q_j)^{\frac{2}{N-2-2a}}<\e
\end{align}
\item[(B)] if $|q_j|>1$
\begin{align}\label{eq:65}
\bigg\|\dfrac{ \tilde u(\tilde u(\tilde
 q_j)^{-\frac{2}{N-2-2a}}x+\tilde q_j)}{\tilde u(\tilde
 q_j)}-z^{a,b}_{K( q_j)}\bigg\|_{\substack{C^{0,\gamma} 
(B_{2R}(0))\cap\\ H^1_a(B_{2R}(0))}}+
|\tilde q_j|\tilde u(\tilde q_j)^{\frac{2}{N-2-2a}}<\e.
\end{align}
\end{itemize}
Moreover the sets
\begin{align*}
U_j:= \begin{cases}
B_{R u(q_j)^{-\frac{2}{N-2-2a}}}(q_j) &\text{in case {\rm (A)}}  \\
\text{Inv}\big(B_{R \tilde u(\tilde q_j)^{-\frac{2}{N-2-2a}}}(\tilde q_j)\big)
&\text{in case {\rm (B)}}
\\ 
\end{cases} 
\text{ are disjoint.} 
\end{align*}
Furthermore, $u$ satisfies
\begin{align*}
u(x) \le C_0 \omega_a(x) \max_{1\le j\le k} \dist(x,q_j)^{-\frac{N-2-2a}{2}}.
\end{align*}
\end{proposition}
\begin{proof}
Fix $\eps>0$ and $R>1$. Let $C_0$ be as in Lemma
\ref{l:standardblowup}. First we apply Lemma~\ref{l:standardblowup} with
${\mathcal K}=\emptyset$ and find $q_1\in \rz^N \cup \{\infty\}$ the
maximum point of  
$u/\omega_a$. If $u(x)\leq C_0 \omega_a(x)\dist
(x,q_1)^{-\frac{N-2-2a}2}$ holds we stop here. Otherwise we apply
again Lemma \ref{l:standardblowup} to obtain $q_2$. From estimates
(\ref{eq:64}) and (\ref{eq:65}) it follows that $U_1$ and $U_2$ are
disjoint. We continue the process. Since $u\in L^p(\R^N,|x|^{-bp})$ and
\begin{align*}
\int\limits_{U_j}\frac{K(x)}{|x|^{bp}}u(x)^p\,dx
&\geq \frac 1{2A_1} \int_{B_R(0)} \frac{\big(z_{K(q_j)}^{a,b}\big)^p}{|y+\e q_j/|q_j||^{bp}}\,dy 
\geq c(a,b,N),
\end{align*}
where $c(a,b,N)$ is independent of $q_j$, $u$, $R>1$ and $\e<1$, we
will stop after a finite number of steps. 
\end{proof}

\begin{proposition}
\label{p:4.2:yyli}
Under the assumptions and notations of Lemma \ref{l:standardblowup} there exist  
for any $0<\e<1$ and $R>1$ some positive constants 
$C_0 =C_0(\e,R,N,a,b,A_1,A_2)$ and $\delta=\delta(\e,R,N,a,b,A_1,A_2)$ such that if $u$ is a solution of
(\ref{eq:49}) with 
\[
\max _{x\in \R^N\cup\{\infty\}}\frac{u(x)}{\omega_a(x)}>C_0 
\]
then 
\[
\dist(q_j,q_\ell) \ge \delta \text{ for all } 1\le j\neq \ell \le k, 
\]
where $q_j=q_j(u)$, $q_\ell=q_\ell(u)$ and $k=k(u)$ are given in Proposition \ref{p:standardblowup}.
\end{proposition}
\begin{proof}
To obtain a contradiction we assume that for some constants $\eps$, $R$, $A_1$ and $A_2$ there exist
sequences ${K_i}$ and $u_i$ satisfying the assumptions of Proposition \ref{p:4.2:yyli} such that
\[ \lim_{i \to \infty}\min_{j \neq \ell}\dist(q_j(u_i),q_\ell(u_i))=0.\]
We may assume that
\begin{align}
\label{eq:28}
\sigma_i:=\dist(q_1(u_i),q_2(u_i))=\min_{j \neq \ell}\dist(q_j(u_i),q_\ell(u_i))\to 0 \text{ as }i \to \infty.  
\end{align}
Let us denote $q_j(u_i)$ by $q_j^i$. Since $U_1(u_i)$ and $U_2(u_i)$
are disjoint and (\ref{eq:28}) holds we have that $u_i(q_1^i) \to 
\infty$ and $u_i(q_2^i) \to \infty$. Therefore we can pass to a subsequence still denoted by $\{u_i\}$
and find $R_i \to \infty$, $\eps_i \to 0$ such that either
$q_1^i=q_1(u_i) \to 0$ or 
$|q_1^i| \to \infty$, and for $j=1,2$ 
\begin{align}
&\bigg\|\frac{u_i\big(u_i(q_j^i)^{\frac2{N-2-2a}}x+q_j^i\big)}{u_i(q_j^i)} 
-z_{K(q_j^i)}^{a,b}\bigg\|_{C^{0,\gamma}(B_{2R_i}(0))}+|q_j^i|u_i(q_j^i)^{\frac2{N-2-2a}}<\e_i
\text{ if }q_1^i\to 0 \label{eq:27}\\
&\bigg\|\frac{\tilde u_i\big(\tilde u_i(\tilde
q_j^i)^{\frac2{N-2-2a}}x+\tilde q_j^i\big)}{\tilde u_i(\tilde q_j^i)}-z_{K(\tilde
q_j^i)}^{a,b}\bigg\|_{C^{0,\gamma}(B_{2R_i}(0))}
+|\tilde q_j^i|\tilde u_i(\tilde q_j^i)^{\frac2{N-2-2a}}<\e_i\ \text{
    if }|q_1^i|\to \infty.\notag
\end{align}
We first
consider the case $q_1^i\to 0$. Since
$U_1(u_i)$ and $U_2(u_i)$ are 
disjoint we have that 
\begin{equation}\label{eq:32}
\sigma_i>c(N) \max_{j=1,2}\{R_i u_i(q_j^i)^{-\frac2{N-2-2a}}\}.
\end{equation}
From (\ref{eq:27}) and (\ref{eq:32}) we get that 
$
\sigma_i^{-1}|q_j^i|<\frac{\e_i}{c(N)R_i}\to 0
$ for $j=1,2$
and obtain the contradiction 
\[
\frac12<|\sigma_i^{-1}(q_2^i-q_1^i)|\to 0.
\]
Performing the same analysis as above for the Kelvin transform $\tilde
u$ of $u$ leads to a
contradiction if $\tilde q_1^i\to 0$. 
\end{proof}

\begin{remark}\label{rem:priori-estimates}
Propositions \ref{p:standardblowup} and \ref{p:4.2:yyli} imply that
there are only finitely many blow-up points and all are isolated.
\end{remark}

\begin{proposition}
\label{p:yyli:thm4.2}
Suppose $\{K_i\}$ and $a\in ](N-4)/2,(N-2)/2[$ satisfy the assumptions of Lemma
\ref{l:standardblowup} and Proposition \ref{p:simplezero}.  
Let $\{u_i\}$ be solutions of (\ref{eq:pi}) with
$\Omega =\rz^N$. Then after passing to a subsequence either $\{u_i/\omega_a\}$
stays bounded in $L^\infty(\rz^N)$ or $\{u_i\}$ has precisely one
blow-up point, which can be at $0$ or at $\infty$.  
\end{proposition}
\begin{proof}
Suppose that $\{u_i/\omega_a\}$ is not uniformly bounded in $L^\infty(\rz^N)$, otherwise
there is nothing to prove. Consequently we may apply Proposition \ref{p:standardblowup} and
Proposition \ref{p:4.2:yyli} to obtain isolated points
$\{q_1^i, \dots, q_{k(i)}^i\}$ satisfying (\ref{eq:64}) and (\ref{eq:65}) with $R_i \to \infty$ and
$\eps_i \to 0$. To obtain a contradiction, we assume that up to a subsequence $k(i)\ge 2$.
Since $u(q_j^i)/\omega_a(q_j^i) \to \infty$ for $j=1,2$ and $\dist(q_1^i,q_2^i)\ge \delta>0$ we may
assume $q_1^i \to 0$ and $q_2^i \to \infty$ and $k(i)=2$ as $i \to \infty$.  
From Proposition \ref{p:simplezero} and Remark
\ref{rem:priori-estimates} they are isolated simple blow-up points. 
From Proposition \ref{p:upperbound} we infer that
\begin{align*}
\lim_{i \to \infty} u_i({q_1^i}) u_i(x) = h(x) \text{ in }
C^{0}_{\rm 
loc}(\rz^N\setminus\{0\}),\\
\div(|x|^{-2a} \nabla h)=0 \text{ in } \rz^N\setminus \{0\}.
\end{align*}
Using Theorem \ref{t:bocher} for $h$ and its Kelvin transform and the
maximum principle we obtain for some $a_1,a_2>0$
\begin{align*}
h(x)= a_1|x|^{2+2a-N} + a_2.  
\end{align*}
We may now proceed as in the proof of Proposition \ref{p:simplezero}
to see that
\begin{align*}
\int_{\partial B_{\sigma}(q_1^i)}B(\sigma,x,u_i,\n u_i) =o(u_i(q_1^i)^{-2}).
\end{align*}
Multiplying by $u_i({q_1^i})^{2}$ and letting $i\to\infty$ we find that
\[
\int_{\partial B_{\sigma}}B(\sigma,x,h,\n h)=0.
\]
This contradicts for small $\sigma$ Proposition \ref{p:boundterm} and completes the proof.
\end{proof}

\begin{proposition}\label{p:deltabound}
Suppose $K\in C^2(\R^N)$ satisfies (\ref{eq:60})-(\ref{eq:51}),
\begin{align*}
a \ge 0,\ \frac{N-4}{2}<a<\frac{N-2}{2}, \text{ and } \frac{4}{N-2-2a}<p<2^*.  
\end{align*}
Then there exists $C_K>0$ such that for any $t\in
(0,1]$ and any solution $u_t$ of
\begin{align}
\label{eq:pt}
\tag{$P_t$}
-\div (|x|^{-2a}\n u)=(1+t\,(K(x)-1))\frac{u^{p-1}}{|x|^{bp}}, \quad u>0\
 \text{in}\ \Di
\end{align}
there holds
\begin{align}\label{eq:8}
\|u_t\|_{E}< C_K 
\end{align}
and 
\begin{align}
\label{eq:9}
C_K^{-1}<u_t\omega_a^{-1}<C_K.  
\end{align}
\end{proposition}
\begin{proof}
The bound in (\ref{eq:9}) follows from (\ref{eq:8}) and Harnack's
inequality in \cite{Gutierrez89}. The estimate in Lemma \ref{l:harnack_potential} of the
appendix shows that $(1+t\,(K(x)-1)) u^{p-2}|x|^{-bp}$ belongs to the required class of
potentials in \cite{Gutierrez89}. 
To show that $u_t/\omega_a$ is bounded in $L^{\infty}(\R^N)$
we argue by contradiction and may  
assume in view of Proposition \ref{p:yyli:thm4.2} that there exists a sequence
$\{t_i\}\subset(0,1]$ converging to $t_0\in [0,1]$ as 
$i\to\infty$ such that $u_{t_i}$ has precisely one blow-up point $(x_i)$,
which can be supposed to be zero using the Kelvin transform.
Corollary \ref{c:pohoRn} yields
\begin{equation*}
0=\int_{\R^N}x\cdot\n K(x)\frac{u_{t_i}(x)^p}{|x|^{bp}}\,dx.
\end{equation*}
Since $0$ is assumed to be the only blow-up point, the Harnack
inequality and (\ref{eq:upperbound}) yield, for any $\sigma\in(0,1)$,
\begin{align*}
\bigg|\int_{B_{\sigma}(x_i)}x\cdot\n
K(x)\frac{u_{t_i}(x)^p}{|x|^{bp}}\,dx\bigg|&=
\bigg|\int_{\R^N\setminus B_{\sigma}(x_i)}x\cdot\n
K(x)\frac{u_{t_i}(x)^p}{|x|^{bp}}\,dx\bigg|\leq C(\sigma)\Big(u_{t_i}(x_i)^{-p}\Big).
\end{align*}
We have that from Taylor expansion, (\ref{eq:16}), and (\ref{eq:7})
\begin{align}\label{eq:41}
|\n K(x_i)|\leq {\rm const} |x_i|=o\Big(u_{t_i}(x_i)^{-\frac2{N-2-2a}}\Big)
\end{align}
and
\begin{align*}
&\bigg|\int_{B_{\sigma}(x_i)}\n K(x)\cdot
x\frac{u_{t_i}^p}{|x|^{bp}}\,dx\bigg|\\
&\quad=\bigg|\int_{B_{\sigma}(x_i)}\n
K(x_i)\cdot
x\frac{u_{t_i}^p}{|x|^{bp}}\,dx+\int_{B_{\sigma}(x_i)}D^2K(x_i)(x-x_i)\cdot
x\frac{u_{t_i}^p}{|x|^{bp}}\,dx\\
&\quad\quad+\int_{B_{\sigma}(x_i)}o(|x-x_i|)\cdot
x\frac{u_{t_i}^p}{|x|^{bp}}\,dx\bigg|.
\end{align*}
From Lemma \ref{l:yy:2.4:zero} and \eqref{eq:41} we infer
 \begin{align*} 
\bigg|\int_{B_{\sigma}(x_i)}\n
K(x_i)\cdot
x\frac{u_{t_i}^p}{|x|^{bp}}\,dx+\int_{B_{\sigma}(x_i)}o(|x-x_i|)\cdot
x\frac{u_{t_i}^p}{|x|^{bp}}\,dx\bigg|=o\Big(u_{t_i}(x_i)^{-\frac4{N-2-2a}}\Big).
\end{align*}
Hence 
\[
\int_{B_{\sigma}(x_i)}D^2K(x_i)(x-x_i)\cdot
x\frac{u_{t_i}^p}{|x|^{bp}}\,dx=o\Big(u_{t_i}(x_i)^{-\frac4{N-2-2a}}\Big).
\]
Since by Lemma \ref{l:yy:2.4:zero} 
\[
\int_{r_i\leq |x-x_i|\leq\sigma}D^2K(x_i)(x-x_i)\cdot
x\frac{u_{t_i}^p}{|x|^{bp}}\,dx=o\Big(u_{t_i}(x_i)^{-\frac4{N-2-2a}}\Big)
\]
we have
\begin{align}\label{eq:42}
\int_{B_{r_i}(x_i)}D^2K(x_i)(x-x_i)\cdot
x\frac{u_{t_i}^p}{|x|^{bp}}\,dx=o\Big(u_{t_i}(x_i)^{-\frac4{N-2-2a}}\Big).
\end{align}
Making in \eqref{eq:42} the change of variables $x=u_{t_i}(x_i)^{-2/(N-2-2a)}y+x_i$ and
using Proposition \ref{prop:Bubble:zero}
\begin{align*}
0=\int_{\R^N}D^2K(0)y\cdot y
{|y|^{-bp}}{z^{a,b}_{1+t_0(K(0)-1)}(y)^p}\,dy=\Delta K(0)\int_{\R^N} {|y|^{2-bp}}{z^{a,b}_{1+t_0(K(0)-1)}(y)^p}\,dy
\end{align*}
which is not possible in view (\ref{eq:7}).
\end{proof}

\begin{pfn}{Theorem \ref{t:compactness}.} It follows from Proposition
  \ref{p:deltabound} and Lemma \ref{l:change}.
\end{pfn}

\noindent
We define $f_{K,\e}:\Di\to\R$ by 
\begin{align*}
f_{K,\e}(u)&=f_0(u)-\e G_K(u)\\
f_0(u)&=\frac12 \int_{\R^N}|x|^{-2a}|\n u|^2-\frac 1p \int_{\R^N}\frac{|u|^{p}}{|x|^{bp}}\\
G_K(u)&=\frac 1p \int_{\R^N}\frac{K(x)|u|^{p}}{|x|^{bp}}.
\end{align*}
We will use the notation $f_{\e}$ (respectively $G$) instead of $f_{K,\e}$ (respectively $G_K$)
whenever there is no possibility of confusion. 
Let us denote by $Z$ the manifold
\[
Z=\{z_{\mu}=z_{1,\mu}^{a,b}:\ \mu>0\}
\]
of the solutions to (\ref{eq:49}) with $K\equiv 1$.
\begin{lemma}\label{p:implicit}
Suppose $p>3$. There exist constants $\rho_0,\e_0,C>0$, and smooth functions 
\begin{align*}
w=w(\mu,\e):\quad &(0,+\infty) \times (-\e_0,\e_0)\ \longrightarrow\
\Di\\
\eta=\eta(\mu,\e):\quad &(0,+\infty) \times (-\e_0,\e_0)\ \longrightarrow\
\R
\end{align*}
such that for any $\mu>0$ and $\e\in(-\e_0,\e_0)$
\begin{align}
\label{eq:12}
w(\mu,\e)\ \text{ is orthogonal to }\ T_{z_{\mu}}Z\\
\label{eq:17}
f_{\e}'\big(z_{\mu}+w(\mu,\e)\big)=\eta(\mu,\e) \dot \xi_{\mu}\\
\label{eq:18}
|\eta(\mu,\e)|+\|w(\mu,\e)\|_{\Di}\leq C\,|\e|\\
\label{eq:26}
\|\dot w (\mu,\e)\|_{\Di}\leq C(1+\mu^{-1})\,|\e|,
\end{align}
where $\dot{\xi}_\mu$ denotes the normalized tangent vector
$\frac{d}{d\mu} z_\mu$ and $\dot{w}$ stands for the derivative of~$w$
with respect to $\mu$.  
Moreover, $(w,\eta)$ is unique in the
sense that there exists $\rho_0>0$
such that if $(v,\tilde\eta)$ satisfies $\|v\|_{\Di}+|\tilde \eta|<\rho_0$
and (\ref{eq:12})-(\ref{eq:17}) for
some $\mu>0$ and $|\eps|\le \eps_0$, then $v=w(\mu,\eps)$ and
$\tilde\eta=\eta(\mu,\eps)$.
\end{lemma}
\begin{proof}
Existence, uniqueness, and estimate \eqref{eq:18} are proved in
\cite{FS1}. In fact 
$w$ and $\eta$ are implicitly defined by
$H(\mu,w,\eta,\eps)=(0,0)$ where 
\begin{align*}
&H: (0,\infty)\times \Di \times \rz
\times \rz \to \Di \times \rz\\
&H(\mu,w,\eta,\eps):= (f_\eps'(z_\mu+w)-\eta \dot{\xi}_\mu,(w,\dot{\xi}_\mu)).
\end{align*}
Let us now show
estimate (\ref{eq:26}). There exists a positive constant $C_*$ such
that for any $\mu>0$ (see \cite{FS1})
\begin{align*}
\bigg\|\bigg(\frac{\rand H}{\rand(w,\eta)}(\mu,0,0,0)\bigg)^{-1}\bigg\| \le
 C_*.
\end{align*}
Since $\dot w$ satisfies
\begin{align*}
\begin{pmatrix}
\dot w\\
\dot \eta
\end{pmatrix}
=-
\bigg(\frac{\partial H}{\partial (w,\eta)}\bigg)^{-1}\bigg|_{(\mu,w,\eta,\eps)}
\cdot\frac{\partial H}{\partial \mu}\bigg|_{(\mu,w,\eta,\eps)}
\end{align*}
we have for $\e$ small using (\ref{eq:18}) and the fact that $f_0 \in C^3$
\begin{align*}
\|\dot w(\mu,\e)\|&\leq C_* \frac{\partial H}{\partial
  \mu}\bigg|_{(\mu,w,\eta,\eps)} \\
&\leq C_* \bigg(\bigg\|f_\eps''(z_\mu+w(\mu,\eps))\dot z_\mu
-\eta(\mu,\eps) \frac{d}{d\mu}\dot\xi_\mu\bigg\| +
\bigg|\bigg(w(\mu,\eps),\frac{d}{d\mu}\dot\xi_\mu\bigg)\bigg|\bigg)\\
&\leq C(1+\mu^{-1})|\eps| + \|f_0''(z_\mu+w(\mu,\eps))\dot z_\mu\|\\
&\le C(1+\mu^{-1})|\eps| + O(\|w(\mu,\eps)\|)\|\dot z_\mu\|\\
&\le C(1+\mu^{-1})|\eps|.
\end{align*}
This ends the proof.
\end{proof}

\begin{corollary}
\label{c:tubo}
Suppose $p>3$ and $K$ satisfies the assumptions of Proposition \ref{p:deltabound}.
Then there exist $t_0>0$ and $R_0>0$ such that any solution $u_t$ of
(\ref{eq:pt}) for 
$t\le t_0$ is of the form $z_\mu + w(\mu,t)$, where $1/R_0<\mu<R_0$.
\end{corollary}
\begin{proof}
First we show that there exists $R_1>0$ and $t_1>0$ such that 
any solution $u_t$ of (\ref{eq:pt}) for $t<t_1$ satisfies
\begin{align*}
\dist(u_t,Z_{R_1})<\rho_0,   
\end{align*}
where by $\dist$ we mean the distance in the $\Di$-norm, $\rho_0$ is
given in Lemma \ref{p:implicit}, and $Z_{R_1}:=\{z_\mu \where 1/R_1
<\mu<R_1\}$. 
By contradiction, assume there exist $R_i\to \infty$, $t_i \to 0$, and solutions $u_{t_i}$ of
(\ref{eq:pt}) such that $\dist(u_{t_i},Z_{R_i})\ge\rho_0$. From (\ref{eq:8}) we can pass to a
subsequence converging weakly in $\Di$ to some $\bar u$; since in view
of the regularity results of \cite{SF_reg}  
$\{u_t\}$ is bounded in $C^{0,\gamma}$ and such a bound excludes any
possibility of concentration, the convergence is actually strong and
$\dist(\bar u,Z)\ge \rho_0$. Furthermore, 
$\bar u$ solves (\ref{eq:pt}) with $t=0$ and hence $\bar u \in Z$, which is impossible.\\
Fix a solution $u_t$ of (\ref{eq:pt}) for some $t< t_1$. A short computation shows
\begin{align*}
\lim_{\mu \to 0}\dist(z_\mu,u_t)^2 = \lim_{\mu \to \infty}\dist(z_\mu, u_t)^2= \|z_1\|^2+\|u_t\|^2>\rho_0^2.  
\end{align*}
Consequently there exists $R_0>0$ independent of $t$ and $z_\mu \in Z_{R_0}$ such that
\begin{align*}
\dist(u_t, Z)= \|u_t-z_\mu\| \text{ and } u_t-z_\mu \in T_{z_\mu}Z^{\perp}.  
\end{align*}
Since $u_t$ solves (\ref{eq:pt}) we have $f'_t(z_\mu+u_t-z_\mu)=0$ and the uniqueness in Lemma
\ref{p:implicit} yields the claim.     
\end{proof}

\section{Leray-Schauder degree}
\noindent We introduce the
Melnikov function  
\begin{equation*}
\Gamma_K(\tau)=\frac 1p \int_{\R^N} K(x) \frac{z_{\tau}^p}{|x|^{bp}}.
\end{equation*}
It is known (for details see \cite{FS1}) that it is possible to extend
the $C^2$- function  
$\Gamma_K$ by continuity to $\tau=0$ and
\begin{align}
\label{eq:43}
\Gamma_K'(0)=0 \text{ and } \Gamma_K''(0)=\frac{\laplace K(0)}{Np}
\int_{\R^N} |x|^2 \frac{z_1(x)^p}{|x|^{bp}}.   
\end{align}
Furthermore, using the Kelvin transform, we find
\begin{align}
\label{eq:44}
\Gamma_K(\tau)=\Gamma_{\tilde K}(\tau^{-1})\quad\text{where}\quad
\tilde K(x)=K(x/|x|^2).  
\end{align}
We define for small $t$ the function $\Phi_{K,t}(\mu):= f_{K,t}(z_{\mu}+w(\mu,t))$ and will denote it
by~$\Phi_t$ whenever there is no possibility of confusion.
\begin{lemma}
\label{l:gamma:critical}
Let $p>3$ and assume $\Gamma_K$ has only non-degenerate critical
points. Then there exists $t_1>0$ 
such that for any $0<t<t_1$ any solution $u_t$ of (\ref{eq:pt}) is of the form
$u_t=z_{\mu_t}+w(\mu_t,t)$, where $\Phi_{K,t}'(\mu_t)=0$ and $\mu_t \in (R_0^{-1},R_0)$ for some
positive $R_0$. Moreover,
up to a subsequence as $t \to 0$
\begin{align}
\label{eq:30}
|\mu_t -\bar\mu|=O(t),  
\end{align}
where $\bar \mu$ is a critical point of $\Gamma_K$. Viceversa, for any critical point $\bar \mu$
of $\Gamma_K$ and for any $0<t<t_1$ there exists one and only one critical point $\mu_t$ of
$\Phi_{K,t}$ such that (\ref{eq:30}) holds.       
\end{lemma}
\begin{proof}
By Corollary \ref{c:tubo} any solution $u_t$ of (\ref{eq:pt}) is of the form
$u_t=z_{\mu_t}+w(\mu_t,t)$, where $\Phi_{t}'(\mu_t)=0$ and $R_0^{-1}< \mu_t<R_0$.
Using the Taylor expansion and (\ref{eq:18}) - (\ref{eq:26}), we have that for  $R_0^{-1}< \mu<R_0$ 
\begin{align}\label{eq:33}
\Phi_{t}'(\mu)&=f_{t}'(z_{\mu}+w(\mu,t))(\dot z_{\mu}+\dot w(\mu,t))\notag\\
              &=f_t'(z_{\mu})(\dot z_{\mu}+\dot w(\mu,t))+(f_t''(z_{\mu})w(\mu,t),\dot z_{\mu}+\dot
              w(\mu,t))+O(\|w(\mu,t)\|^2))\notag\\
              &=-t G'(z_{\mu})(\dot z_{\mu}+\dot w(\mu,t))+(f_0''(z_{\mu})w(\mu,t),\dot
              w(\mu,t)))\notag\\
              &\quad-t (G''(z_{\mu})w(\mu,t),\dot z_{\mu}+\dot
              w(\mu,t))+O(\|w(\mu,t)\|^2))\notag\\
              &=-t\Gamma'(\mu)+O(t^2).
\end{align}
Fix a sequence $(t_n)$ converging to $0$. Since $\mu_t$ is bounded, we may assume that
$(\mu_{t_n})$ converges to $\bar \mu$. From expansion (\ref{eq:33}) we have that
\[
0=\Phi_{t_n}'(\mu_{t_n})=-t_n(\Gamma'(\mu_{t_n})+O(t_n))
\]
hence $\bar \mu$ is a critical point of $\Gamma$. A further expansion yields
\[
0=\Phi_{t_n}'(\mu_{t_n})-t_n(\Gamma''(\bar \mu)(\mu_{t_n}-\bar\mu)+o(\mu_{t_n}-\bar\mu))+O(t_n^2)
\]
which gives for $n\to\infty$
\[
(\mu_{t_n}-\bar\mu)(\Gamma''(\bar\mu)+o(1))=O(t_n)
\]
proving (\ref{eq:30}) for $\Gamma''(\bar\mu)\not= 0$. Viceversa let $\bar \mu$ be a critical point
of $\Gamma$. Arguing as above we find as $\mu\to\bar\mu$ and for any $0<t<t_1$
\[
\Phi_{t}'(\mu)=t
(\mu-\bar\mu)(\Gamma''(\bar\mu)+o(1))+O(t^2)
\]
hence there exists $\mu_t$ such that
\[
\mu_t=\bar \mu-(\Gamma''(\bar\mu)+o(1))^{-1}O(t)\quad\text{and}\quad \Phi_{t}'(\mu_t)=0.
\]
To prove uniqueness of such a $\mu_t$, we follow \cite{Bapo} and expand $\Phi_t$ in a critical point
$\mu_t$
\begin{align}\label{eq:34}
\Phi_t''(\mu_t)&=\big(f_t''(z_{\mu_t}+w(\mu_t,t))(\dot z_{\mu_t}+\dot
w(\mu_t,t)),(\dot z_{\mu_t}+\dot 
w(\mu_t,t))\big) \notag\\
& =\big(f_0''(z_{\mu_t}+w(\mu_t,t))(\dot z_{\mu_t}+\dot
w(\mu_t,t)),(\dot z_{\mu_t}+\dot 
w(\mu_t,t))\big) \notag\\
&\quad -t \big(G''(z_{\mu_t}+w(\mu_t,t))(\dot z_{\mu_t}+\dot
w(\mu_t,t)),(\dot z_{\mu_t}+\dot 
w(\mu_t,t))\big) \notag\\
&=(f_0''(z_{\mu_t})\dot w(\mu_t,t),\dot w(\mu_t,t))+(f_0'''(z_{\mu_t})
w(\mu_t,t)(\dot z_{\mu_t}+\dot 
w(\mu_t,t)),\dot 
z_{\mu_t} +\dot w(\mu_t,t) \notag)\\
&\quad -t \big(G''(z_{\mu_t}+w(\mu_t,t))(\dot z_{\mu_t}+\dot
w(\mu_t,t)),(\dot z_{\mu_t}+\dot 
w(\mu_t,t))\big) \notag\\
&= (f_0'''(z_{\mu_t}) w(\mu_t,t) \dot z_{\mu_t},\dot z_{\mu_t})
-t \big(G''(z_{\mu_t})\dot z_{\mu_t},\dot z_{\mu_t}\big) + O(t^2).
\end{align}
Since any critical point $\mu_t$ of $\Phi_t$ gives rise to a critical
point $z_{\mu_t}+w(\mu_t,t)$ of $f_t$, we have that
\begin{align}\label{eq:35}
0&=(f_t'(z_{\mu_t}+w(\mu_t,t)),\ddot z_{\mu_t})\notag\\
 &=(f_t'(z_{\mu_t})+f_t''(z_{\mu_t}) w(\mu_t,t)+O(\|
 w(\mu_t,t)\|^2),\ddot z_{\mu_t})\notag\\ 
 &=-t (G'(z_{\mu_t}),\ddot z_{\mu_t})+(f_0''(z_{\mu_t})w(\mu_t,t),\ddot
 z_{\mu_t}) +O(t^2).
\end{align}
Differentiating $f_0''(z_{\mu_t})\dot z_{\mu_t}=0$ and testing with
$w(\mu_t,t)$ we obtain 
\begin{align}\label{eq:36}
0=(f_0'''(z_{\mu_t})\dot z_{\mu_t}\dot
z_{\mu_t},w(\mu_t,t))+(f_0''(z_{\mu_t})\ddot z_{\mu_t} ,w(\mu_t,t)).
\end{align}
Putting together (\ref{eq:35}) and (\ref{eq:36}) we get
\begin{align*}
(f_0'''(z_{\mu_t})\dot z_{\mu_t}\dot
z_{\mu_t},w(\mu_t,t))=-t (G'(z_{\mu_t}),\ddot z_{\mu_t})+O(t^2)
\end{align*}
 hence in view of (\ref{eq:34}) 
\begin{align}\label{eq:37}
\Phi_t''(\mu_t)= -t (G'(z_{\mu_t}),\ddot z_{\mu_t})
-t \big(G''(z_{\mu_t})\dot z_{\mu_t},\dot z_{\mu_t}\big) +
O(t^2)=-t\Gamma''(\mu_t)+O(t^2). 
\end{align}
To prove uniqueness, we choose $\delta>0$ such that $\sgn
\Gamma''(\mu)=\sgn \Gamma''(\bar\mu)\not=0$ for any
$|\mu-\bar\mu|<\delta$. From (\ref{eq:37}), there exists $t(\delta)>0$
such that if 
$t<t(\delta)$ and $\mu_t$ is a critical point of $\Phi_t$ such that 
 $|\mu_t-\bar\mu|<\delta$, then 
\[
\sgn \Phi_t''(\mu_t)=-\sgn \Gamma''(\bar\mu).
\]
From (\ref{eq:33}) we have that  for $t<t(\delta)$
\begin{align*}
\sgn\Gamma''(\bar\mu)&=\deg(\Gamma',B_{\delta}(\bar\mu),0)=
\deg(-\Phi_t',B_{\delta}(\bar\mu),0)\\
&=-\sum_{\substack{y\in
    B_{\delta}(\bar\mu)\\   \Phi_t'(y)=0}}\sgn\Phi_t''(y)=\#\{y\in
B_{\delta}(\bar\mu):\  \Phi_t'(y)=0\}\,\sgn\Gamma''(\bar\mu).
\end{align*}
Hence $\#\{y\in
B_{\delta}(\bar\mu):\  \Phi_t'(y)=0\}=1$, proving uniqueness. 
\end{proof}
\noindent
\begin{lemma}\label{l:id-comp}
For any $K\in L^{\infty}(\R^N)$ the operator 
\begin{align*}
L_K:\ u\mapsto \big(-\div
(|x|^{-2a}\n)\big)^{-1}\frac{K(x)}{|x|^{bp}}|u|^{p-2}u 
\end{align*}
is compact from $E$ to $E$.
\end{lemma}
\begin{proof}
Let $\{u_n\}$ be a bounded sequence in $E$ and set $v_n=L_K(u_n)$, i.e.
\[
-\div (|x|^{-2a}\n v_n)=\frac{K(x)}{|x|^{bp}}|u_n|^{p-2}u_n.
\]
By Caffarelli-Kohn-Nirenberg inequality, $\{v_n\}$ is bounded in
$\Di$ and passing to a subsequence we may assume that it converges
weakly in $\Di$ and pointwise almost everywhere to some limit
$v\in\Di$. Since $\{u_n\}$ is uniformly bounded in
$L^{\infty}(B_3(0))$, from \cite{SF_reg} the sequence $\{v_n\}$ is
uniformly bounded in $C^{0,\gamma}(B_2(0))$. Using the Kelvin
transform we arrive at
\begin{align*}
-\div (|x|^{-2a}\n \tilde v_n)&=|x|^{-(N+2+2a)+bp}K(x/|x|^2)
\,|u_n(x/|x|^2)|^{p-2}u_n(x/|x|^2)\\ 
& = K(x/|x|^2)\frac{|\tilde u_n|^{p-2}\tilde u_n}{|x|^{bp}}.
\end{align*}
Since $\{u_n\}$ is uniformly bounded in $E$,
$\{\tilde u_n\}$ is uniformly bounded in 
$L^{\infty}(B_3(0))$ and hence from \cite{SF_reg} the sequence
$\{\tilde v_n\}$ is
uniformly bounded in $C^{0,\gamma}(B_2(0))$. Since a uniform bound in
$C^{0,\gamma}(B_2(0))$ implies equicontinuity and
\[
\|(v_n-v_m)\omega_a^{-1}\|_{C^{0}(\R^N\setminus B_1(0))}\leq {\rm
  const}\,\|\tilde v_n-\tilde v_m\|_{C^0(B_1(0))}
\]
from the
Ascoli-Arzel\`a Theorem there exists a subsequence
$\{v_{n}\}$ strongly converging in $C^0(\R^N,\omega_a)$ to
$v$. Moreover, the $C^0(\R^N,\omega_a)$-convergence excludes any
possibility of concentration at $0$ or at $\infty$ and $\{v_n\}$
converges strongly in $\Di$.
\end{proof}
\noindent 
From Proposition \ref{p:deltabound}, there exists a positive constant
$C_K$ such that $\|u\|_E< C_K$ and $C_K^{-1}<u\omega_a^{-1}$
 for any solution $u$ of
(\ref{eq:pt}) uniformly with respect to $t\in (0,1]$. By the
above lemma, the Leray-Schauder degree $\deg(Id-L_K,{\mathcal B}_K,0)$
 is well-defined, where ${\mathcal B}_K:=\{u\in E:\,\|u\|_E<C_K,
C_K^{-1}<u\omega_a^{-1}\}$. 
\begin{theorem}\label{t:grado} Under the assumptions of Proposition
  \ref{p:deltabound} and for $p>3$ we have
\[
\deg(Id-L_K,{\mathcal B}_K,0)=-\frac{\sgn\Delta K(0)+\sgn 
  \Delta \tilde K(0)}2.
\]
\end{theorem}
\begin{proof}
By transversality, we can assume that $\Gamma_K$ has only
non-degenerate critical points. If not, we proceed with a small
perturbation of $K$. By Proposition \ref{p:deltabound} and the
homotopy invariance of the Leray-Schauder degree, for $0<t<t_1$ 
\begin{align*}
\deg(Id-L_K,{\mathcal B}_K,0)=\deg(Id-L_{tK},{\mathcal B}_K,0).
\end{align*}
By Lemma \ref{l:gamma:critical} we have
\begin{align*}
\deg(Id-L_{tK},{\mathcal
  B}_K,0)&=\sum_{\mu\in(\Phi_{t,K}')^{-1}(0)}(-1)^{{\mathfrak
    m}(z_{\mu}+w(\mu,t),f_{t,K})}
\end{align*}
where ${\mathfrak
    m}(z_{\mu}+w(\mu,t),f_{t,K})$ denotes the Morse index of $f_{t,K}$
  in $z_{\mu}+w(\mu,t)$. 
We will only sketch the computation of ${\mathfrak m}(z_{\mu}+w(\mu,t),f_{t,K})$ and refer
to \cite{AmBa1,Bapo,kromer} for details.    
The spectrum of $f_0''(z_\mu)$ is completely known (see \cite{FS1}) and $\Di$ is decomposed
in $\langle z_\mu \rangle \oplus T_{z_\mu}Z \oplus \langle z_\mu,T_{z_\mu}Z \rangle^\perp$, where $z_\mu$ is an
eigenfunction of $f_0''(z_\mu)$ with corresponding eigenvalue $-(p-2)$,
$T_{z_\mu}Z=\ker(f_0''(z_\mu))$, and $f_0''(z_\mu)$ restricted to the orthogonal
complement of $\langle z_\mu,T_{z_\mu}Z \rangle$ is bounded below by a positive constant. 
Consequently, to compute the Morse index ${\mathfrak m}(z_{\mu}+w(\mu,t),f_{t,K})$ for
small $t$ it is enough to know the behavior of $f''_{t,K}(z_{\mu}+w(\mu,t))$ along
$T_{z_\mu}Z$. From the expansion
\begin{align*}
f_{t,K}(z_{\mu}+w(\mu,t)) = f_0(z_\mu) - t \Gamma_K(\mu) + o(t^2)= {\rm const}\,- t
\Gamma_K(\mu) + o(t^2)    
\end{align*}
we have that for $t$ small  
\begin{equation}\label{eq:38}
{\mathfrak m}(z_{\mu}+w(\mu,t),f_{t,K})=1+
\begin{cases}
1&\text{if }\Gamma_K''(\mu)>0\\
0&\text{if }\Gamma_K''(\mu)<0.
\end{cases}
\end{equation}
From \eqref{eq:38} and Lemma \ref{l:gamma:critical}, we know that for
$t$ small  
\begin{align*}
\sum_{\mu\in(\Phi_{t,K}')^{-1}(0)}(-1)^{{\mathfrak
    m}(z_{\mu}+w(\mu,t),f_{t,K})}&=-\sum_{\mu\in(\Gamma_{K}')^{-1}(0)}
(-1)^{{\mathfrak m}(\mu,-\Gamma_K)}\\
&=\deg(\Gamma_K',((R_0+1)^{-1},R_0+1),0),  
\end{align*}
where $R_0$ is given in Lemma \ref{l:gamma:critical}. From (\ref{eq:43}) we obtain for $\mu \to 0$ 
\begin{align*}
\Gamma_K'(\mu)= \Gamma_K''(0)\mu + o(\mu) =  {\rm const}\laplace K(0)\mu + o(\mu). 
\end{align*}
Hence $\sgn \Gamma_K'((R_0+1)^{-1})= \sgn\laplace K(0)$. Using (\ref{eq:44}) for obtain for $\mu \to
\infty$
\begin{align*}
\Gamma_K'(\mu)= -\mu^{-2}\Gamma'_{\tilde K}(\mu^{-1})=  -{\rm const}\laplace \tilde K(0)\mu^{-3}
+o(\mu^{-3}). 
\end{align*}
Therefore $\sgn \Gamma_K'((R_0+1))= -\sgn\laplace \tilde K(0)$ and 
\begin{align*}
\deg(\Gamma_K',((R_0+1)^{-1},R_0+1),0)= -\frac{\sgn\Delta K(0)+\sgn 
  \Delta \tilde K(0)}2,  
\end{align*}
which proves the claim.
\end{proof}

\begin{pfn}{Theorem \ref{t:existence}.} It follows directly from
  Theorem \ref{t:grado} and Lemma \ref{l:change}.
\end{pfn}

\appendix

\section{}

\begin{lemma}\label{l:change}
$v$ is a solution to (\ref{eq:45}) if and only if
$u(x)=|x|^{a-\a}v(x)$ solves (\ref{eq:49}), where $a=a(\a,\l)$ and $b=b(\a,\b,\l)$ are
given in \eqref{eq:ab}.
\end{lemma}
\begin{proof}
By standard elliptic regularity $u$ and $v$ are $C^2(\rz^N\setminus\{0\})$.
Consequently we may compute for $x \in \rz^N\setminus\{0\}$
\begin{align*}
\div(|x|^{-2a}\n
  u(x))=(a-\a)(N-a-\a-2)|x|^{-a-\a-2}v(x)+|x|^{\a-a}\div(|x|^{-2\a}\n
  v)
 \end{align*}
and hence, in view of (\ref{eq:45})
\begin{align*}
-\div(|x|^{-2a}\n 
  u(x))
&=\big[\l+(\a-a)(N-2-\a-a)\big]\frac{u(x)}{|x|^{2a+2}}+K(x)\frac{u^{p-1}}
{|x|^{p(a-\a+\b)}}
.
\end{align*}
From \eqref{eq:ab} we have that $\l+(\a-a)(N-2-\a-a)=0$ and
$a-\a+\b=b$. 
Since $C^\infty(\rz^N\setminus\{0\})$ is dense in ${\mathcal{D}^{1,2}_{\a}(\R^N)}$ and 
$\Di$ (see \cite{CatrinaWang}), the lemma is thereby proved.
\end{proof}
\begin{lemma}\label{l:decrease}
Let $\{K_i\}_i$ satisfy (\ref{eq:blow1}), $(u_i)_{i \in \nz}$
satisfy (\ref{eq:pi}) and 
$x_i \to 0$ be an isolated blow up point. Then for any $R_i \to
\infty$, there exists a positive constant $C$
depending on $\lim_{i \to \infty}K_i(x_i)$ and $a$, $b$, and $N$ such
that  after passing to a subsequence the function
  $\bar w_i$ defined in 
  (\ref{eq:62}) is strictly decreasing for $C
  u_i(x_i)^{-2/(N-2-2a)}\le r\le r_i$ where $r_i=R_i u_i(x_i)^{-\frac
    2{N-2-2a}}$. 
\end{lemma}
\begin{proof}
Making the change of variable $y=u_i(x_i)^{-\frac 2{N-2-2a}}x+x_i$,
there results
\begin{align*}
\bar w_i(r)&=
\frac{r^{\frac{N-2-2a}2}}{|\partial
  B_r(x_i)|}\int_{\partial B_r(x_i)}\ui(y)\\
&=r^{\frac{N-2-2a}2}\media_{\partial B_{r
    u_i(x_i)^{ 2/(N-2-2a)}}(0)}\ui\big(u_i(x_i)^{-\frac
  2{N-2-2a}}x+x_i\big).
\end{align*}
From the proof of Proposition \ref{prop:Bubble:zero} we
    have that for some function $g_i\in C^{0,\gamma}(B_{2R_i}(0))$
$$\ui\big(u_i(x_i)^{-\frac
  2{N-2-2a}}x+x_i\big)=\ui(x_i)\big(z_{K(0)}^{a,b}(x)+g_i(x)\big)
$$
where $\|g_i\|_{C^2(B_{2R_i}(0)\setminus
  B_{C}(0))}\leq \e_i$.
Being $z_{K(0)}^{a,b}$ a radial function, from above we find
\begin{align*}
\bar w_i(r)&=r^{\frac{N-2-2a}2}\ui(x_i)\media_{\partial B_{r
    u_i(x_i)^{ 2/(N-2-2a)}}(0)}\big(z_{K(0)}^{a,b}(x)+g_i(x)\big)\\
&=r^{\frac{N-2-2a}2}\ui(x_i)\big[z_{K(0)}^{a,b}(r
    u_i(x_i)^{ 2/(N-2-2a)})+{\textstyle\mediapiccola_{\partial B_{r
    u_i(x_i)^{ 2/(N-2-2a)}}(0)}g_i}\big].
\end{align*}
A direct computation shows that
\begin{align*}
&\frac{d}{d r}\bar w_i(r)\\
&\ =\ui(x_i)
    r^{\frac{N-4-2a}2}(z_{K(0)}^{a,b}(r
    u_i(x_i)^{ 2/(N-2-2a)}))^{\frac
      p2}\Big[{\textstyle\frac{N-2-2a}2}\Big(1-K(0)\ui(x_i)^{p-2} 
      r^{\frac{(p-2)(N-2-2a)}2}\Big)\\
&\quad +{\textstyle\frac{N-2-2a}2 }\big({\textstyle\mediapiccola
  g_i}\big)z_{K(0)}^{a,b}(r 
    u_i(x_i)^{ 2/(N-2-2a)})^{-\frac p2}+r({\textstyle\mediapiccola
      g_i})'z_{K(0)}^{a,b}(r 
    u_i(x_i)^{ 2/(N-2-2a)})^{-\frac p2}
\Big].
\end{align*}
Since for $C u_i(x_i)^{-2/(N-2-2a)}\le r\le r_i$, there results $C\leq r
    u_i(x_i)^{ 2/(N-2-2a)}\leq R_i$, we have that 
\[
\media_{\partial B_{r
    u_i(x_i)^{ 2/(N-2-2a)}}(0)}g_i\leq \e_i,\quad\frac{d}{d r}\media_{\partial B_{r
    u_i(x_i)^{ 2/(N-2-2a)}}(0)}g_i\leq \e_i.
\]
Moreover for
$C=\big(\frac{1+\delta}{K(0)}\big)^{\frac2{(p-2)(N-2-2a)}}$ we have
$1-K(0)\ui(x_i)^{p-2} 
      r^{\frac{(p-2)(N-2-2a)}2}\leq -\delta$. Choosing
      $\e_i=o\Big(R_i^{-\frac {p(N-2-2a)}2}\Big)$ the claim follows.
\end{proof}

\begin{lemma}
\label{l:harnack_potential}
Suppose $a,b,p$ satisfy (\ref{eq:47}) and (\ref{eq:ab}). Let $(z_i)_{i\in \nz} \subset
\rz^N$ and consider the measures $\mu_i:=|x-z_i|^{-2a} \,dx$, then we have  
for $0<r<2$ as $r \to 0$
\begin{align*}
\sup_{x \in B_2(0),\, i \in \nz} \int_{B_r(x)} |y-z_i|^{-bp} \int_{|x-y|}^8 \frac{s\, ds}{\mu_i(B_s(x))}\, dy
\to 0.  
\end{align*}
\end{lemma}
\begin{proof}
We use as $c$ a generic constant that may change its value from line to line.
Fix $x \in B_2(0)$. From the doubling property of the measure $\mu_i$ (see \cite{potential}) we find
\begin{align*}
M_i(x,|x-y|) &:= \int_{|x-y|}^8 \frac{s\, ds}{\mu_i(B_s(x))}\\
&\le c\, \begin{cases}
|x-y|^{-N+2a+2}, &\text{if }|x-y|>\frac{1}{2}|x-z_i|\\
|x-y|^{-N+2} |x-z_i|^{2a}+|x-z_i|^{-N+2a+2}, &\text{if }|x-y|\le\frac{1}{2}|x-z_i|.
\end{cases}
\end{align*}
An easy calculation shows that $2a-bp>-2$ and that if $a\ge 0$ then $2a-bp\le 0$. 
Hence, we may estimate for $0<r \le \frac{1}{2}|x-z_i|$ and $y \in B_r(x)$
\begin{align*}
|y-z_i|& \ge |x-z_i|-|x-y| \ge \frac{1}{2}|x-z_i|
\end{align*}
and 
\begin{align*}
\int_{B_r(x)} |y-z_i|^{-bp} M_i(x,|x-y|) \,dy &\le c r^{2+2a-bp}.     
\end{align*}
Since $-bp>-2-2a>-N$ we may use the above estimate to derive
\begin{align*}
\intg_{B_{2|x-z_i|}(x)} &|y-z_i|^{-bp} M_i(x,|x-y|) \,dy \le c |x-z_i|^{2+2a-bp}.
\end{align*} 
Consequently we obtain for $\frac{1}{2}|x-z_i|\le r \le 2|x-z_i|$
\begin{align*}
\int_{B_{r}(x)} |y-z_i|^{-bp} M_i(x,|x-y|) \,dy
&\le c |x-z_i|^{2+2a-bp} \le c r^{2+2a-bp}.
\end{align*}
Finally we obtain for $2|x-z_i|< r \le 2$ and $|x-y|>2|x-z_i|$
\begin{align*}
|y-z_i|& \ge |y-x|-|x-z_i| \ge \frac{1}{2}|y-x|
\end{align*}
and 
\begin{align*}
\int_{B_r(x)} &|y-z_i|^{-bp} M_i(x,|x-y|) \,dy 
\le c r^{2+2a-bp},     
\end{align*}
which ends the proof.  
\end{proof}

\noindent A function $u$ will be called $\mu$-harmonic in $\Omega\subset \rz^N$, if $u \in
D_{a,loc}^{1,2}(\rz^N)$ and for all $\varphi \in C^\infty_c(\Omega)$
there holds
\begin{align*}
\int_\Omega |x|^{-2a} \nabla u \nabla \varphi =0.  
\end{align*}
Let us prove a B\^ocher-type theorem for $\mu$-harmonic functions.
\begin{theorem}
\label{t:bocher} Let $u$ be a nonnegative $\mu$-harmonic function in
$\R^N\setminus\{0\}$. 
Then there exist a constant $A\geq0$ and a H\"older continuous
function $B$, $\mu$-harmonic in $\R^N$, such that
\[
u(x)=A\,|x|^{2+2a-N}+B(x).
\]
\end{theorem}
\begin{proof}
We distinguish two cases.\par
\noindent
{\bf Case 1:} there exists a sequence $x_n\to0$ and a positive
constant $M$ such that
$|u(x_n)|\leq M$. In this case 
the Harnack Inequality (Theorem 6.2 of \cite{potential}) implies that
$u$ is bounded. Moreover from \cite[Lemma 6.15]{potential} $u$ can be
continuously extended to $0$ and is a weak solution of
\[
-\div(|x|^{-2a}\n u)=0\quad\text{in }\R^N,
\]
see \cite[Lemma 2.1]{CatrinaWang}. Therefore from the Liouville
Theorem 
\cite [Theorem 6.10]{potential} we get that $u$ is constant and the
theorem holds 
with  $A=0$ and $B\equiv {\rm const}$.
\par
\noindent
{\bf Case 2:} $u(x_n)\to+\infty$ for any sequence $x_n\to0$.
We can extend $u$ in $0$ to be $u(0):=+\infty$, thus obtaining a lower
semi-continuous function in $\R^N$. Moreover
\cite [Theorem 7.35]{potential} implies that $u$ is super-harmonic in
the sense of the 
definition of \cite[Chapter 7]{potential}, i.e. 
\begin{itemize}
\item[(i)] $u$ is lower semi-continuous,
\item[(ii)] $u\not\equiv\infty$ in each component of $\R^N$,
\item[(iii)] for each open $D\Subset \R^N$ and each $h\in
  C^0(\R^N)$ $\mu$-harmonic in $D$ the inequality $u\geq h$ on
  $\partial D$ implies $u\geq h$ in
  $D$.
\end{itemize}
Let us remark that in order to apply Theorem 7.35 in
\cite{potential} we need to prove that $0$ has capacity $0$ with
respect to our weight; indeed
\begin{align*}
{\rm cap}_{|x|^{-2a}}(\{0\},\R^N)&:=\inf_{\substack {u\in
  C^{\infty}_0(\R^N),\ u\equiv1\\
  \text{in a neighborhood of}\ 0}}\int_{\R^N}|x|^{-2a}|\n u|^2\leq 
{\rm cap}_{|x|^{-2a}}(B_r,\R^N)\\
&\leq{\rm
  cap}_{|x|^{-2a}}(B_r,B_{2r})\leq c r^{N-2-2a}
\end{align*}
for any $r>0$, where we have used \cite[Lemma
2.14]{potential}. Then ${\rm
  cap}_{|x|^{-2a}}(\{0\},\R^N)=0$. From 
\cite[Corollary 7.21]{potential} there holds 
\[
-\div(|x|^{-2a}\n u)\geq0\quad\text{in the sense of distributions on}\ \R^N
\]
hence from the Riesz Theorem there exists a Radon measure $\mu\geq0$
in $\R^N$ such that 
\[
\langle -\div(|x|^{-2a}\n u), \varphi \rangle
=\int_{\R^N}\varphi \,d\mu\quad\forall\, \varphi\in C^{\infty}_0(\R^N).
\]
Since $\langle -\div(|x|^{-2a}\n u), \varphi \rangle=0$ for any 
$\varphi\in C^{\infty}_0(\R^N\setminus\{0\})$,
$\mu$ must be supported in $\{0\}$ and so $\mu=A\d_0$ for a
nonnegative constant
$A$. Since the Green's function $G_a(x):=|x|^{2+2a-N}$ satisfies
\[
-\div(|x|^{-2a}\n G_a)=\d_0\quad\text{in the sense of distributions on}\ \R^N,
\]
we have that
\[
-\div(|x|^{-2a}\n (u-A G_a))=0
\]
in the sense of distributions on $\R^N$. Theorem 3.70 and Lemma 6.47 in
\cite{potential} imply that $B:=u-AG_a$ is H\"older continuous. 
\end{proof}

\bibliographystyle{adinat}
\bibliography{blowup}

\def\cprime{$'$}
\begin{thebibliography}{18}
\expandafter\ifx\csname natexlab\endcsname\relax\def\natexlab#1{#1}\fi
\expandafter\ifx\csname url\endcsname\relax
  \def\url#1{{\tt #1}}\fi
\expandafter\ifx\csname urlprefix\endcsname\relax\def\urlprefix{URL }\fi

\bibitem[{Abdellaoui et~al.(2003)Abdellaoui, Felli and Peral}]{AbdelFelliPeral}
B.~Abdellaoui, V.~Felli and I.~Peral.
\newblock {\em Existence and multiplicity for perturbations of an equation
  involving {H}ardy inequality and critical {S}obolev exponent in the whole
  $\bf{R}\it\sp {N}$\/} (2003).
\newblock Preprint.

\bibitem[{Ambrosetti and Badiale(1998{\natexlab{a}})}]{AmBa2}
A.~Ambrosetti and M.~Badiale.
\newblock {\em Homoclinics: {P}oincar\'e-{M}elnikov type results via a
  variational approach\/}.
\newblock Ann. Inst. H. Poincar\'e Anal. Non Lin\'eaire, {\bf 15}
  (1998{\natexlab{a}}), no.~2, 233--252.

\bibitem[{Ambrosetti and Badiale(1998{\natexlab{b}})}]{AmBa1}
A.~Ambrosetti and M.~Badiale.
\newblock {\em Variational perturbative methods and bifurcation of bound states
  from the essential spectrum\/}.
\newblock Proc. Roy. Soc. Edinburgh Sect. A, {\bf 128} (1998{\natexlab{b}}),
  no.~6, 1131--1161.

\bibitem[{Badiale and Pomponio(2002)}]{Bapo}
M.~Badiale and A.~Pomponio.
\newblock {\em Bifurcation results for semilinear elliptic problems in
  $\bf{R}\it\sp {N}$\/} (2002).
\newblock Preprint.

\bibitem[{Caffarelli et~al.(1984)Caffarelli, Kohn and Nirenberg}]{CKN}
L.~Caffarelli, R.~Kohn and L.~Nirenberg.
\newblock {\em First order interpolation inequalities with weights\/}.
\newblock Compositio Math., {\bf 53} (1984), no.~3, 259--275.

\bibitem[{Catrina and Wang(2001{\natexlab{a}})}]{CatrinaWang}
F.~Catrina and Z.-Q. Wang.
\newblock {\em On the {C}affarelli-{K}ohn-{N}irenberg inequalities: sharp
  constants, existence (and nonexistence), and symmetry of extremal
  functions\/}.
\newblock Comm. Pure Appl. Math., {\bf 54} (2001{\natexlab{a}}), no.~2,
  229--258.

\bibitem[{Catrina and Wang(2001{\natexlab{b}})}]{CatrWang01}
F.~Catrina and Z.-Q. Wang.
\newblock {\em Positive bound states having prescribed symmetry for a class of
  nonlinear elliptic equations in {$\bold R\sp N$}\/}.
\newblock Ann. Inst. H. Poincar\'e Anal. Non Lin\'eaire, {\bf 18}
  (2001{\natexlab{b}}), no.~2, 157--178.

\bibitem[{Chou and Chu(1993)}]{ChouChu}
K.~S. Chou and C.~W. Chu.
\newblock {\em {On the best constant for a weighted Sobolev-Hardy
  inequality}\/}.
\newblock J. Lond. Math. Soc., II. Ser., {\bf 48} (1993), no.~1, 137--151.

\bibitem[{Felli and Schneider(2003{\natexlab{a}})}]{FS1}
V.~Felli and M.~Schneider.
\newblock {\em Perturbation results of critical elliptic equations of
  {C}affarelli-{K}ohn-{N}irenberg type\/}.
\newblock Journal of Differential Equations, {\bf 191} (2003{\natexlab{a}}),
  121--142.

\bibitem[{Felli and Schneider(2003{\natexlab{b}})}]{SF_reg}
V.~Felli and M.~Schneider.
\newblock {\em Regularity results for elliptic equations of
  {C}affarelli-{K}ohn-{N}irenberg type\/} (2003{\natexlab{b}}).
\newblock Preprint.

\bibitem[{Guti{\'e}rrez(1989)}]{Gutierrez89}
C.~E. Guti{\'e}rrez.
\newblock {\em Harnack's inequality for degenerate {S}chr\"odinger
  operators\/}.
\newblock Trans. Amer. Math. Soc., {\bf 312} (1989), no.~1, 403--419.

\bibitem[{Heinonen et~al.(1993)Heinonen, Kilpel{\"a}inen and
  Martio}]{potential}
J.~Heinonen, T.~Kilpel{\"a}inen and O.~Martio.
\newblock {\em Nonlinear potential theory of degenerate elliptic equations\/}.
\newblock The Clarendon Press Oxford University Press, New York (1993).

\bibitem[{Kr\"omer(2002)}]{kromer}
S.~Kr\"omer.
\newblock {\em Symmetriebrechung bei {V}ariationsproblemen\/} (2002).
\newblock Diplomarbeit, Universit\"at Augsburg.

\bibitem[{Li(1995)}]{li95}
Y.~Li.
\newblock {\em {Prescribing scalar curvature on ${\mathbb S}\sp n$ and related
  problems. I.}\/}.
\newblock J. Differ. Equations, {\bf 120} (1995), no.~2, 319--410.

\bibitem[{Schoen and Zhang(1996)}]{SchoenZhang}
R.~Schoen and D.~Zhang.
\newblock {\em Prescribed scalar curvature on the {$n$}-sphere\/}.
\newblock Calc. Var. Partial Differential Equations, {\bf 4} (1996), no.~1,
  1--25.

\bibitem[{Schoen(1991)}]{schoen91}
R.~M. Schoen.
\newblock {\em {On the number of constant scalar curvature metrics in a
  conformal class}\/}.
\newblock In {\em {Differential geometry. A symposium in honour of Manfredo do
  Carmo, Proc. Int. Conf., Rio de Janeiro/Bras. 1988, Pitman Monogr. Surv. Pure
  Appl. Math. 52, 311-320}\/} (1991) .

\bibitem[{Smets(2001)}]{smets}
D.~Smets.
\newblock {\em Nonlinear {S}chroedinger equations with hardy type potential and
  critical nonlinearities\/} (2001).
\newblock Preprint.

\bibitem[{Terracini(1996)}]{terracini}
S.~Terracini.
\newblock {\em On positive entire solutions to a class of equations with a
  singular coefficient and critical exponent\/}.
\newblock Adv. Differential Equations, {\bf 1} (1996), no.~2, 241--264.

\end{thebibliography}

\end{document}